%
%
%
%
%

\newcount\chno  
\newcount\equno 
\newcount\refno 

\font\chhdsize=cmbx12 at 14.4pt

%
\def\startbib{\def\biblio{\bigskip\medskip
  \noindent{\chhdsize References.}\medskip\bgroup\parindent=2em}}
\def\endbib{\edef\biblio{\biblio\egroup}}
\def\reflbl#1#2{\global\advance\refno by 1
  \edef#1{\number\refno}
    \global\edef\biblio{\biblio\par\item{[\number\refno]}#2\par}}

\def\eqlbl#1{\global\advance\equno by 1
  \global\edef#1{{\number\chno.\number\equno}}
  (\number\chno.\number\equno)}

\def\qed{\hfil\hbox to 0pt{}\ \hbox to 2em{\hss}\ 
         \hbox to 0pt{}\hskip-2em plus 1fill
         \vrule height3pt depth1pt width3pt\par\medskip}
\def\eqed{\hfil\hbox to 0pt{}\ \hbox to 2em{\hss}\ 
          \hbox to 0pt{}\hskip-2em plus 1fill
          \vbox{\hrule height .25pt depth 0pt width 7pt
            \hbox{\vrule height 6.5pt depth 0pt width .25pt
              \hskip 6.5pt\vrule height 6.5pt depth 0pt width .25pt}
            \hrule height .25pt depth 0pt width 7pt}\par\medskip}


\def\msimp#1#2{#1%
  \hbox to 0pt{\hskip 0pt minus 3fill
    \phantom{$#1$}\hbox to 0pt{\hss$#2$\hss}\phantom{$#1$}%
    \hskip 0pt minus 1fill}}


\def\pr{{\prime}}
\def\ppr{{\prime\prime}}

\def\({\left(}
\def\){\right)}
\def\[{\left[}
\def\]{\right]}
\def\<{\left\langle}
\def\>{\right\rangle}

\def\ol{\overline}

\def\wtilde{\widetilde}

\def\eps{\epsilon}

\def\sech{{\rm \, sech\,}}

\font\tendouble=msbm10 
\font\sevendouble=msbm7  
\font\fivedouble=msbm5

\newfam\dbfam
\textfont\dbfam=\tendouble \scriptfont\dbfam=\sevendouble
\scriptscriptfont\dbfam=\fivedouble

\mathchardef\dbA="7041 
\mathchardef\dbB="7042 
\mathchardef\dbC="7043 
\mathchardef\dbD="7044 
\mathchardef\dbE="7045 
\mathchardef\dbF="7046 
\mathchardef\dbG="7047 
\mathchardef\dbH="7048 
\mathchardef\dbI="7049 
\mathchardef\dbJ="704A 
\mathchardef\dbK="704B 
\mathchardef\dbL="704C 
\mathchardef\dbM="704D 
\mathchardef\dbN="704E 
\mathchardef\dbO="704F 
\mathchardef\dbP="7050 
\mathchardef\dbQ="7051 
\mathchardef\dbR="7052 \def\RR{{\fam=\dbfam\dbR}}
\mathchardef\dbS="7053 \def\SS{{\fam=\dbfam\dbS}}
\mathchardef\dbT="7054 
\mathchardef\dbU="7055 
\mathchardef\dbV="7056 
\mathchardef\dbW="7057 
\mathchardef\dbX="7058 
\mathchardef\dbY="7059 
\mathchardef\dbZ="705A

\def\cA{{\cal A}}
\def\cK{{\cal K}}
\def\cX{{\cal X}}

\def\rhs{{right-hand side}}
\def\lhs{{left-hand side}}

\def\supp{{\rm supp\,}}
\def\dist{{\rm dist\,}}

\def\arcosh{\rm arcosh\,}

\def\ts{\tilde{s}}
\def\dd{\,{\rm d}}

\def\dr{\dd r}
\def\ds{\dd s}
\def\dts{\dd \tilde{s}}
\def\dss{\dd s^\pr}
\def\dt{\dd t}
\def\dx{\dd x}
\def\dy{\dd y}

\def\dtau{\dd\tau}

\def\dsigma{\dd\sigma}


\magnification= 1200 
\baselineskip= 12pt

\input epsf 


\startbib
\refno=0


\reflbl\Ahlfors{
Ahlfors, L.V., 
	{\it An extension of Schwarz's Lemma},
		Trans. Amer. Math. Soc. {\bf 43}, pp. 359--364 (1938).}


\reflbl\Aubin{
Aubin, T., 
	{\it Meilleures constantes dans le th\'eo\-r\`eme d'inclusion
	de Sobolev et un th\'eo\-r\`eme de Fredholm non lin\'eaire
        pour la transformation conforme de la courbure scalaire},
        	J. Functional Anal. {\bf 32}, pp. 148--174 (1979).}

\reflbl\Aviles{
Aviles, P.,
	{\it Conformal complete metrics with prescribed non-negative
	Gaussian curvature in $\RR^2$},
		Invent. Math. {\bf 83}, pp 519--544 (1986).}

\reflbl\Bandle{
Bandle, C.,
	{\it Isoperimetric Inequalities and Applications,}
		Pitman, Boston, (1980).}

\reflbl\BrezisMerle{ 
Brezis, H., and  Merle, F., 
	{\it Uniform estimates and blow-up behavior of solutions of 
	$-\Delta u = V(x) e^u$ in two dimensions,}
		Commun.\ PDE {\bf 16}, pp. 1223--1253 (1991).}



\reflbl\ChaKieCMP{
Chanillo, S., and Kiessling, M. K.-H.,
	{\it Rotational symmetry of solutions of some nonlinear 
	problems in statistical mechanics and geometry}, 
		Commun. Math. Phys. {\bf 160}, 217--238 (1994). }

\reflbl\ChaKieGAFA{
Chanillo, S., and Kiessling, M. K.-H.,
	{\it Conformally invariant systems of nonlinear 
	PDE of Liouville type},
		Geom. Functional Anal., {\bf 5}, pp. 924--947 (1995).}

\reflbl\ChaKieDUKE{
Chanillo, S., and Kiessling, M. K.-H.,
	{\it Surfaces with prescribed Gauss curvature},
		Duke Math. J. {\bf 105}, pp. 309--353 (2000).}

\reflbl\ChaLi{
Chanillo, S., and Li, Y. Y., 
       	{\it Continuity of solutions of uniformly 
        elliptic equations in $\RR^2$},
	        Manuscr. Math. {\bf 77}, pp.  415--433 (1992).}

\reflbl\ChenLiA{
Chen, W., and  Li, C., 
	{\it Classification of solutions of some nonlinear 
	elliptic equations,}
		Duke Math.\ J. {\bf 63},  pp. 615--622 (1991).} 

\reflbl\ChenLiB{
Chen, W., and Li, C.,
	{\it Qualitative Properties of solutions to some 
	nonlinear elliptic equations in $\RR^2$,}
		Duke Math.\ J. {\bf 71},  pp. 427--439 (1993).} 

\reflbl\ChengLinMA{
Cheng, K.-S., and Lin, C.-S.,
	{\it On the asymptotic behavior of solutions of the 
	conformal Gaussian curvature equations on $\RR^2$}, 
		Math. Ann. {\bf 308}, pp. 119--139 (1997).}

\reflbl\ChengLinJDE{
Cheng, K.-S., and Lin, C.-S.,
	{\it On the conformal Gaussian curvature equations in $\RR^2$}, 
		J. Diff. Eq. {\bf 146}, pp. 226--250 (1998).}

\reflbl\ChengLinASNSP{
Cheng, K.-S., and Lin, C.-S.,
	{\it Compactness of conformal metrics with positive Gaussian 
	curvature in $\RR^2$}, 
		Ann. Scuola Norm. Sup. Pisa Cl. Sci. {\bf 26}, 
		pp. 31--45 (1998).}

\reflbl\ChengLinNA{
Cheng, K.-S., and Lin, C.-S.,
	{\it Conformal metrics in $\RR^2$ with prescribed
	Gaussian curvature with positive total curvature},
		Nonl. Anal. {\bf 38}, pp. 775-783 (1999).}

\reflbl\ChengNiI{
Cheng, K.-S., and Ni, W.-M.,
	{\it On the structure of the conformal Gaussian 
	curvature equation on $\RR^2$}, 
		Duke Math.\ J. {\bf 62},  pp. 721--737 (1991).} 

\reflbl\ChengNiII{
Cheng, K.-S., and Ni, W.-M.,
	{\it On the structure of the conformal Gaussian 
	curvature equation on $\RR^2$, II}, 
		Math. Ann. {\bf 290},  pp. 671--680 (1991).} 

\reflbl\ChipotShafrirWolansky{
Chipot, M., Shafrir I., and Wolansky, G.,
	{\it On the solutions of Liouville systems},
	J. Diff. Eq. {\bf 140}, pp. 59--105 (1997).}

\reflbl\ChouWan{
Chou, K.S., and Wan, T.Y.H.,
	{\it Asymptotic radial symmetry for
	solutions of $\Delta u+e^u=0$ in a punctured disc,}
		Pac. J. Math. {\bf 163}, pp.269--276 (1994).}



\reflbl\GidasNiNirbrg{ 
Gidas, B., Ni, W.-M., and Nirenberg, L.,
	{\it Symmetry and related properties via the maximum principle,}
		Commun.\ Math.\ Phys. {\bf 68}, pp. 209--243 (1979).}

\reflbl\GilbargTrudinger{
Gilbarg, D. and  Trudinger, N.S., 
	{\it Elliptic Partial Differential Equations of Second Order,}
		 Springer Verlag, New York (1983).}



\reflbl\JostWang{
Jost, J., and Wang, G., 
	{\it Analytic aspects of the Toda system: I. 
		A Moser-Trudinger inequality}, 
	e-print http://xxx.lanl.gov/list/math-ph/0011039}



\reflbl\KiePHYSICA{
Kiessling, M.K.-H., 
	{\it Statistical mechanics approach to some problems in
	conformal geometry},  
		Physica A {\bf 79}, pp. 353--368 (2000).}

\reflbl\KieJllPHPL{
Kiessling, M.K.-H., and Lebowitz, J.L.,
	{\it Dissipative stationary plasmas: Kinetic modeling, 
	 Bennett's pinch and  generalizations,}
		The Physics of Plasmas {\bf 1}, pp. 1841--1849 (1994).}

\reflbl\Liouville{
Liouville, J., 
	{\it Sur l'\'equation aux diff\'erences partielles
        $\partial^2 \log \lambda/\partial u\partial v \pm \lambda/ 2 a^2 =0$},
	        J. de Math. Pures Appl. {\bf 18}, pp. 71--72 (1853).}





\reflbl\McOwen{
McOwen, R.C.,
	{\it Conformal metrics in $\RR^2$ with prescribed Gaussian
	curvature and positive total curvature},
		Indiana Univ. Math. J. {\bf 34}, pp. 97--104 (1985).}

\reflbl\Ni{
Ni, W.M., 
	{\it On the elliptic equation $\Delta u + Ke^{2u} =0$ and
	conformal metrics with prescribed Gaussian curvatures},
		Invent. Math. {\bf 66}, pp. 343--352 (1982).}



\reflbl\Osserman{
Osserman, R., 
	{\it On the inequality $\Delta u \geq f(u)$},
		Pac. J. Math. {\bf 7} pp. 1641--1647 (1957).}

\reflbl\Poincare{
Poincar\'e, H., 
       	{\it Les fonctions Fuchsiennes et l'\'equation $\Delta u = e^u$},
	        Journal de Math\'e\-ma\-tiques Pures et Applic\'ees,
		Ser. 5, {\bf 4}, pp. 137--230 (1898).}

\reflbl\PrajapatTarantello{
Prajapat, J., and Tarantello, G.,
	{\it On a class of elliptic problems in $\RR^2$: symmetry and
	uniqueness results},
		Proc. Roy. Soc. Edinburg (to appear).}

\reflbl\Sattinger{
Sattinger, D.H.,
	{\it Conformal metrics in $\RR^2$ with prescribed curvatures},
		Indiana Univ. Math. J. {\bf 22}, pp. 1--4 (1972).}



\reflbl\Tarantello{
Tarantello, G.,
	{\it Vortex condensation of a non-relativistic Chern-Simons theory},
		J. Diff. Eq. {\bf 141}, pp. 295--309 (1997).}

\reflbl\Wittich{
Wittich, H., 
	{\it Ganze L\"osungen der Differentialgleichung $\Delta u=e^u$},
		Math. Z. {\bf 49}, pp. 579--582  (1944).}


\endbib
\centerline{\chhdsize THE CONFORMAL PLATE BUCKLING EQUATION}
\bigskip

\centerline{SAGUN CHANILLO and MICHAEL KIESSLING}
\smallskip
\centerline{Department of Mathematics, Rutgers University}
\centerline{110 Frelinghuysen Rd., Piscataway, NJ 08854}
\bigskip
\bigskip
\bigskip
\bigskip
\bigskip

\noindent
{\bf ABSTRACT}: 
	The linear equation $\Delta^2 u = 1$ for the infinitesimal
buckling under uniform unit load of a thin elastic plate over $\RR^2$ 
has the particularly interesting nonlinear generalization 
${\Delta_g}^2 u = 1$, where $\Delta_g=e^{-2u}\Delta$ is
the Laplace--Beltrami operator for the metric $g = e^{2u}g_0$, 
with $g_0$ the standard Euclidean metric on $\RR^2$. 
	This conformal elliptic PDE of fourth order is equivalent to the 
nonlinear system of elliptic PDEs of second-order
$\Delta u(x)+K_g(x)\exp(2u(x))=0$ 
and 
$\Delta K_g(x) + \exp(2u(x))=0$, 
with  $x\in\RR^2$, describing a conformally flat surface with a
Gauss curvature function $K_g$ that is generated self-consistently 
through the metric's conformal factor.
	We study this conformal plate buckling equation
under the hypotheses of finite integral curvature 
$\int K_g\exp(2u)\dx = \kappa$, 
finite  area 
$\int \exp(2u)\dx = \alpha$,
and the mild compactness condition 
$K_+\in L^1(B_1(y))$, uniformly w.r.t. $y\in\RR^2$.
	We show that asymptotically for $|x|\to\infty$ all 
solutions behave like 
$u(x) = -(\kappa/2\pi)\ln |x| + C +o(1)$ 
and 
$K(x) = - (\alpha/2\pi) \ln |x| + C + o(1)$, 
with  
$\kappa \in (2\pi, 4\pi)$ and $\alpha  = \sqrt{2\kappa(4\pi - \kappa)}$.
	We also show that for each $\kappa \in (2\pi, 4\pi)$ there exists 
a $K^*$ and a radially symmetric solution pair $u,K$, satisfying
${\cal K}(u)=\kappa$ and max$K = K^*$, which is unique modulo translation 
of the origin, and scaling of $x$ coupled with a translation of $u$.
\hfill

\bigskip\bigskip

\centerline{Version of April 17, 2001.}
\vfill

\hrule\medskip
\noindent
$\msimp{\copyright}{c}$ (2001) The authors. 
Reproduction for non-commercial purposes of this article, 
in its entirety and by any means, is permitted.

\vfill\eject

\chno=1
\equno=0
\noindent
{\bf I. INTRODUCTION}
\smallskip

\noindent
In this paper we study the nonlinear, fourth-order elliptic PDE
$$
	{\Delta_g}^2 u(x) = \lambda ;\qquad x\in \RR^2
\eqno\eqlbl\PDE
$$
for a smooth scalar function $u: \RR^2\to \RR$, where 
$\Delta_g = e^{-2u}\Delta_{g_0}$ is the Laplace-Beltrami operator 
w.r.t. the conformally flat metric $g = e^{2u}g_0$, with $g_0$ the
Euclidean standard metric of $\RR^2$ and $\Delta_{g_0}\equiv \Delta$ 
the standard Laplacian on $\RR^2$, and $\lambda \in \RR^+$ a parameter.
	In the limit of small $u$, the nonlinear equation (\PDE)
reduces to the linear equation
$$
	{\Delta}^2 u(x) = \lambda;\qquad x\in \RR^2\, ,
\eqno\eqlbl\PDEplate
$$
which is familiar from the linear theory of stationary buckling 
of a thin, elastic plate under uniform load $\lambda$.
	For this reason, we will call (\PDE) 
{\it the conformal plate buckling equation}.

	For fixed $\lambda$, equation (\PDE)
is invariant under the isometries of Euclidean space $\RR^2$ and
under the scaling $x\mapsto k x$   combined with the translation 
$u\mapsto u -\ln k$, where $k>0$.
	On the punctured plane (\PDE) is invariant also under the 
Kelvin transform (inversion) $x \mapsto x/|x|^2$ combined with the map
$u(x) \mapsto u(x/|x|^2) - 2\ln |x|$. 
	However, as we shall see, the singularity at the origin is not
removable so that invariance under the  full Euclidean group of
$\RR^2$ does not hold. 

	If we allow $\lambda$ to change its value under a
transformation, then (\PDE) is invariant also under the combined
transformation $u\mapsto u +u_0$, and $\lambda\mapsto e^{-4u_0}\lambda$.
	Thus, by choosing the constant $u_0 = \ln\lambda^{1/4}$ we can 
always achieve that
$$
\lambda = 1.
\eqno\eqlbl\gauge
$$
	Henceforth we assume  (\gauge) without loss of generality.

	For $\lambda =1$ the fourth-order equation (\PDE) 
is equivalent to the nonlinear system of second-order elliptic PDEs
$$
\eqalignno{
-\Delta u(x) & = K(x)e^{2u(x)}\, ,
&\eqlbl\PDEu\cr
-\Delta K(x) & =  e^{2u(x)}\, ,
&\eqlbl\PDEK}
$$
which describes a conformally flat surface over $\RR^2$
with metric $g = e^{2u}g_0$ and Gauss curvature function $K\equiv K_g$ 
generated in a self-consistent manner.
	While a considerably literature has accumulated 
about the celebrated prescribed Gauss curvature problem
where $K$ is given and only $u$ has to be found by solving (\PDEu), 
see [\Aubin, \Aviles, \BrezisMerle, \ChaKieCMP, \ChaKieDUKE, 
\ChenLiA, \ChenLiB, \ChengLinMA, \ChengLinJDE, \ChengLinASNSP, 
\ChengLinNA, \ChengNiI, \ChengNiII, \KiePHYSICA, \McOwen, \Ni,
\PrajapatTarantello, \Sattinger] and further references therein,
the literature on self-consistent Gauss curvature problems is relatively 
sparse [\ChaKieGAFA, \ChipotShafrirWolansky, \JostWang, \KieJllPHPL].
	In particular, we are not aware of any previous study of the 
self-consistent Gauss curvature problem (\PDEu), (\PDEK),
equivalently  the conformal plate buckling equation.

	We now present our main results for the conformal plate
buckling equation, which we state in their equivalent self-consistent
Gauss curvature form.
	We are interested in an infinite surface
with  finite  area
$$
\cA(u) = \int_{\RR^2}e^{2u(x)}\dx\, 
\eqno\eqlbl\AREA
$$
and finite integral curvature
$$
	\cK(u) 
= 
	\int_{\RR^2 }K(x)e^{2u(x)}{\rm d}x\,.
\eqno\eqlbl\iCURV
$$
\vfill\eject

\smallskip

\noindent
{\bf Theorem 1.1:} 	
{\it Assume $u\in C^{2,\alpha}$ and $K \in C^{2,\alpha}$ jointly solve 
(\PDEu) and (\PDEK) for finite integral curvature, $\cK(u)=\kappa$, 
and finite  area, $\cA(u)=\alpha$.
	In addition assume that $K_+ \in L^1(B_1(x_0))$ uniformly
w.r.t. $x_0$, where $K_+\equiv \max\{K,0\}$.
	Then, uniformly as $|x|\to\infty$, we have
$$
\eqalignno{
	u(x) 
&= 
	- \kappa {1\over 2\pi} \ln |x| + u(0) +
	{1\over 2\pi} \int_{\RR^2}\ln|y| K(y) e^{2u(y)}\dy +o(1)\, ,
&\eqlbl\lapuint\cr
	K(x) 
&=
	-\alpha {1\over 2\pi} \ln |x| + K(0) +
{1\over 2\pi} \int_{\RR^2} \ln|y| e^{2u(y)}\dy +o(1)\, ,
 &\eqlbl\lapKint}
$$
with $\kappa \in (2\pi,4\pi)$, and with $\alpha\in (0,2^{3/2}\pi)$ given by
$$
\alpha =\sqrt{2\kappa(4\pi - \kappa )}\, .
\eqno\eqlbl\pohozaev
$$
}
\smallskip

\noindent
{\it Remarks}: 
	1. Since $\kappa\in (2\pi,4\pi)$, the map $\kappa\mapsto\alpha$ 
given in (\pohozaev) is strictly monotonic decreasing, hence invertible, 
so that alternately to (\pohozaev) we have
$$
	\kappa 
= 
	2\pi\left(1 
	+ \sqrt{1-{1\over 2}\left({\alpha\over 2\pi}\right)^2}\right). 
\eqno\eqlbl\pohoALT
$$ 
	2. The corresponding results for general positive load $\lambda$
in (\PDE) obtain by replacing $\alpha\mapsto \sqrt{\lambda}\,\alpha$ 
in (\pohozaev) and (\pohoALT). 
	This leaves the bounds on $\kappa$ unchanged, i.e.
$2\pi< \kappa<4\pi$, while the bounds on $\alpha$ change to
$0<\alpha < \sqrt{2/\lambda}\, 2\pi$. 
\eqed

	Our next theorem asserts that the range of integral curvature
values $\kappa\in (2\pi,4\pi)$ displayed in  Theorem 1.1  is optimal,
and so is then the associated range of values of the  area 
$\alpha\in (0,2^{3/2}\pi)$.
\smallskip

\noindent
{\bf Theorem 1.2:} {\it For each $\kappa\in (2\pi,4\pi)$ there exists
a value $K^* >0$ and 
a  pair of $C^\infty$ functions $u,K$ which is 
radially symmetric and decreasing about some point $x_*$, which jointly 
solves (\PDEu), (\PDEK), and for which $\cK(u) = \kappa$ and $K(x_*)=K^*$.
	This solution pair is unique up to translations of $x_*$, and 
scalings $x\mapsto kx$ coupled with the translations $u\mapsto u-\ln k$.}
\smallskip

\noindent
{\it Remark:} A typical solution pair $u,K$ is illustrated in 
3 figures at the end of the paper.\eqed

	We conclude our introduction with two interesting open questions.
\smallskip

\noindent
{\bf Open Problem 1.3:} 
{\it Is the value $K^*$ in Theorem 1.2 uniquely determined by each
$\kappa\in (2\pi,4\pi)$?}
\smallskip

	We can show that there is a surjective map $K^*\mapsto \kappa$
on the interval of admissible $K^*$; Open Problem 1.3 asks whether this 
map is also injective.
\smallskip

\noindent
{\bf Open Problem 1.4:} 
{\it Given the conditions stated in Theorem 1.1, are all
     solutions $u,K$ of (\PDEu), (\PDEK) radially symmetric?} 
\smallskip

	We tend to believe that the answer to Open Problem 1.4
is affirmative, but so far a proof has resisted all our attempts.
\smallskip

	We now turn to the proofs of our two theorems. 
	Theorem 1.1 will be proved in section 2 essentially by
harmonic analysis techniques. 
	Theorem 1.2 is proved in section 3 by mapping the ODE's 
for the radial solutions to a scattering problem of a Newtonian 
point particle in $\RR^2$ and applying techniques from potential 
scattering theory.

\vfill\eject

\medskip
\chno=2
\equno=0
\noindent
{\bf II. PROOF OF THEOREM 1.1.}
\smallskip

\noindent
	We begin with the observation that standard 
elliptic theory tells us that, if $u$ and $K$ jointly solve
(\PDEu) and (\PDEK), with $u\in C^{2,\alpha}$, then 
by (\PDEK) also $K\in C^{4,\alpha}$, from 
which it now follows via (\PDEu) that $u\in C^{4,\alpha}$, whence 
$u\in C^{\infty}$ and $K\in C^{\infty}$ by bootstrapping.  

We next state a representation lemma.
\smallskip

\noindent
{\bf Lemma 2.1}: 
{\it Together with the hypotheses of Theorem 1.1, equations
(\PDEu) and (\PDEK) are equivalent to the pair of integral equations
$$
\eqalignno{
	u(x) 
 &= 
	u(0) -{1\over 2\pi}	\int_{\RR^2} 
		\Big(\ln |x-y| -\ln|y|\Big)K(y)e^{2u(y)}
			\dy ,
&\eqlbl\inturep\cr
	K(x)  
&= 
	K(0) -{1\over 2\pi}	\int_{\RR^2} 
		\Big(\ln |x-y| -\ln|y|\Big) e^{2u(y)}
			\dy\, .
&\eqlbl\intKrep}
$$
}
\medskip

{\it Proof of Lemma 2.1}: 
	Clearly, if $u,K$ jointly solve (\inturep), (\intKrep)
and satisfy the other hypotheses of Theorem 1.1, then
$u,K$ jointly solve (\PDEu), (\PDEK) under these hypotheses. 
	To prove the converse, let $u\in C^{\infty}$ satisfy 
$\int \exp(2u) \dx < \infty$, and let $K\in C^{\infty}$ solve (\PDEK). 
	Then $K$ is given by
$$
	K(x)  
= 
	H(x) -{1\over 2\pi}\int_{\RR^2} 
		\Big(\ln |x-y| -\ln|y|\Big) e^{2u(y)}\dy\, ,
\eqno\eqlbl\GREENSintK
$$
where $H(x)$ is an entire harmonic function on $\RR^2$. 
	Now, by hypothesis, $K_+\in L^1(B_1(x_0))$, uniformly 
w.r.t. $x_0\in\RR^2$.
	Thus, from (\GREENSintK) and $\exp(2u)\in L^1(\RR^2)$, 
we have that $H(x) \leq C + C\ln |x|$, whence $H$ is a constant.
	By inspection of (\GREENSintK) it now follows that $H = K(0)$.  
	
	We now take into account that our $u$ also solves (\PDEu), 
and that $\int K \exp(2u) \dx < \infty$. 
	Then $u$ is given by
$$
	u(x)  
= 
	h(x) -{1\over 2\pi}	\int_{\RR^2} 
		\Big(\ln |x-y| -\ln|y|\Big) K(y) e^{2u(y)}\dy\, ,
\eqno\eqlbl\GREENSintu
$$
where $h(x)$ is an entire harmonic function on $\RR^2$. 
	We now show that $h(x) = u(0)$.

	To this effect, having just proved (\intKrep), we now 
observe that (\intKrep) tells us that $K(x)<0$ for $|x| > R$ 
(with $R$ sufficiently large, depending on $u$), whence $u$ is  
subharmonic for $|x|>R$, and so is $u_+$, the positive part of $u$.
	Thus, for $|y|>2R$ we have 
$\| u_+\|_{L^\infty(B_{1/2}(y))} \leq C \| u_+\|_{L^1(B_1(y))}$, 
with $C$ independent of $y$ for $|y|>2R$.
	But then, since $u\in C^\infty$, we even have 
$\| u_+\|_{L^\infty(B_{1/2}(y))} \leq C \| u_+\|_{L^1(B_1(y))}$, 
with $C$ independent of $y\in\RR^2$.
	Furthermore, we have $\| u_+\|_{L^1(B_1(y))}<C$ uniformly w.r.t.
$y\in \RR^2$.
	Namely, setting $\Lambda_y =\supp u_+ \cap B_1(y)$, we have
$\| u_+\|_{L^1(B_1(y))}= \| u\|_{L^1(\Lambda_y)} 
\leq \int_{\Lambda_y}\exp(2u)\dx 
\leq \int_{\RR^2} \exp(2u)\dx < \infty$, the last step by our hypothesis.
	Thus, we conclude that $\| u_+\|_{L^1(B_1(y))}< C$ 
uniformly w.r.t. $y\in \RR^2$, as claimed.
	Hence, $\| u_+\|_{L^\infty(B_{1/2}(y))} \leq C$ 
uniformly w.r.t. $y\in \RR^2$, i.e. $ u_+ \in L^\infty(\RR^2)$.
	Finally, from $ u_+ \in L^\infty(\RR^2)$, together with
(\GREENSintu) and $K\exp(2u)\in L^1(\RR^2)$, we conclude that 
$h(x) \leq C + C\ln |x|$, whence $h$ is a constant, $h = u(0)$
by inspection of (\GREENSintu). \qed
\medskip

\noindent
{\bf Corollary 2.2}: 
{\it Assume $u,K$ jointly solve  (\PDEu), (\PDEK)  
and satisfy the other hypotheses of Theorem 1.1.
	Then, uniformly as $|x|\to \infty$, we have 
$$
\eqalignno{
	u(x) 
&= 
	- \kappa {1\over 2\pi} \ln |x| + u(0) +
	{1\over 2\pi} \int_{\RR^2}\ln|y| K(y) e^{2u(y)}\dy +o(1),
&\eqlbl\uASYMP\cr
	K(x) 
&=
	- \alpha {1\over 2\pi}\ln |x| + K(0) +
{1\over 2\pi} \int_{\RR^2} \ln|y| e^{2u(y)}\dy +o(1)\, . 
&\eqlbl\Kasymp}
$$
}

\smallskip
{\it Proof of Corollary 2.2:} 
	By Lemma 2.1, $u,K$ jointly solve (\inturep), (\intKrep),
with $\cK(u) =\kappa$ and $\cA(u) =\alpha$.
	By (\intKrep) and $\cA(u) =\alpha$, we immediately have
$$
\lim_{|x|\to\infty} {K(x)\over \ln|x|} = - {1\over 2\pi}\alpha.
\eqno\eqlbl\Kln
$$
	Since furthermore $\cK(u) =\kappa$, we now conclude that
$\int_{\RR^2}\ln(1+|x|)\exp(2u(x))\dx < \infty$. 
	With these estimates our Corollary 2.2 now follows at once
from (\inturep), (\intKrep). \qed
\medskip

\noindent
{\bf Corollary 2.3}: {\it Under the hypotheses of Theorem 1.1, the
integral curvature is bounded below by 
$$
\kappa > 2\pi.
\eqno\eqlbl\kBOUNDlow
$$
}

\medskip
{\it Proof of Corollary 2.3}: 
	Assume $\kappa \leq 2\pi$. 
	It then follows immediately from the asymptotic formula
(\uASYMP) that  $\int_{\RR^2}\exp(2u)\dx = \infty$, in contradiction
to our hypothesis that $\cA(u) =\alpha$. 
	Hence, the lower bound (\kBOUNDlow) follows.\qed
\medskip

	Our next result is a Pokhozaev identity for the system 
(\PDEu), (\PDEK).
\smallskip

\noindent
{\bf Proposition 2.4}: 
{\it Under the hypotheses of Theorem 1.1, the integral curvature $\kappa$
and the  area $\alpha$ 
satisfy the identity
$$
	\alpha^2 
= 
	2\kappa( 4\pi - \kappa).
\eqno\eqlbl\POHOuKC
$$
}
\smallskip

{\it Proof of Proposition 2.4}:
	We multiply (\PDEu) by $-x\cdot\nabla u(x)$ and
(\PDEK) by $- x\cdot\nabla K(x)$, then integrate over $B_R$,
apply the usual scheme of integrations by parts on the left-hand sides, 
and get, respectively, 
$$
	R\int_{\partial B_R}\!
		 \left(\big(\nu \cdot \nabla u(x)\big)^2
		- {1\over 2}|\nabla u(x)|^2\right)\!
	\dsigma 
= 
	- {1\over 2} \int_{B_R} K(x) x\cdot \nabla e^{2u(x)} \dx,
\eqno\eqlbl\POHOuA
$$
$$
	R\int_{\partial B_R}\!
	\left(
		\big(\nu \cdot \nabla K(x)\big)^2
		- {1\over 2}|\nabla K(x)|^2
	\right)\!
	\dsigma 
= 
	-\int_{B_R}e^{2u(x)}  x\cdot \nabla K(x)\dx .
\eqno\eqlbl\POHOKA
$$
	By multiplying (\POHOKA) by $1/2$ and adding
the result to (\POHOuA) we obtain 
$$
\eqalignno{
	R\int_{\partial B_R}\!
	\left(
		\big(\nu \cdot \nabla u(x)\big)^2
		- {1\over 2}|\nabla u(x)|^2
		+ {1\over 2}\big(\nu \cdot \nabla K(x)\big)^2
		- {1\over 4}|\nabla K(x)|^2
	\right)
	\!\dsigma ,
\qquad &
&\cr
=	-{1\over 2} \int_{B_R} 
		x\cdot \nabla \left(K(x)e^{2u(x)}\right)
	&\dx . \qquad
&\eqlbl\POHOuKA}
$$
	Integrating next by parts on the right-hand side, using that
$\nabla\cdot x =2$ for $x\in \RR^2$, and moving
the resulting surface integral over to the left-side, we get
$$
\eqalignno{
&	R\int_{\partial B_R}\!
	\left(
		\big(\nu \cdot \nabla u(x)\big)^2
		- {1\over 2}|\nabla u(x)|^2
		+ {1\over 2}\big(\nu \cdot \nabla K(x)\big)^2
		- {1\over 4}|\nabla K(x)|^2
		+ K(x)e^{2u(x)}
	\right)\!
	 \dsigma 
&\cr
&\qquad\qquad\qquad\qquad\qquad\qquad\qquad\qquad\qquad\qquad\qquad\qquad
=	
	\int_{B_R} K(x)e^{2u(x)}\dx .
&\eqlbl\POHOuKB}
$$
	We now let $R\to\infty$.
	Clearly,
$$
	\int_{B_R} K(x)e^{2u(x)} \dx
\to \kappa\quad	
	\qquad{\rm as}\quad R\to \infty.
$$
	Furthermore, from Corollary 2.2 we infer right away that
$$
\eqalignno{
	R\int_{\partial B_R} K(x)e^{2u(x)} \dsigma 
\to 0\quad	
	\qquad{\rm as}\quad R\to &\infty
\, , &\eqlbl\asyONE\cr
	R\int_{\partial B_R}\!
	\left(
		\big(\nu \cdot \nabla u(x)\big)^2
		- {1\over 2}|\nabla u(x)|^2
	\right)\!
	 \dsigma 
\to 
	{\kappa^2\over 4\pi}
	\qquad{\rm as}\quad R\to &\infty
\, , &\eqlbl\asyTWO\cr
	R\int_{\partial B_R}\!
	\left(
		 \big(\nu \cdot \nabla K(x)\big)^2
		- {1\over 2}|\nabla K(x)|^2
	\right)\!
	 \dsigma 
\to 
	{\alpha^2\over 4\pi}
	\qquad{\rm as}\quad R\to &\infty\, .
&\eqlbl\asyTRE}
$$
	Thus, taking the limit $R\to \infty$ in our identity (\POHOuKB)
we obtain (\POHOuKC).
	Since $\alpha >0$, we see that (\POHOuKC) is identical to
(\pohozaev).\qed
\smallskip

\noindent
{\bf Corollary 2.5}: 
{\it The integral curvature is bounded above by
$$
\kappa <4\pi.
\eqno\eqlbl\kBOUNDup
$$ 
	The  area is bounded above by 
$$
\alpha < 2^{3/2}\pi.
\eqno\eqlbl\lBOUNDup
$$}
\smallskip

{\it Proof of Corollary 2.5}:
	The bound (\kBOUNDup) immediately spins off (\POHOuKC),
recalling that, by definition, $\alpha >0$. 
	The bound (\lBOUNDup) is  an immediate consequence of (\POHOuKC) 
and the lower bound $\kappa > 2\pi$, see (\kBOUNDlow) in Corollary 2.3.\qed
\smallskip

This concludes the proof of Theorem 1.1.\qed
\vfill\eject

\medskip
\chno=3
\equno=0
\noindent
{\bf III. PROOF OF THEOREM 1.2.}
\medskip
	In this section we prove the existence of radial solutions $u,K$ 
of the system (\PDEu), (\PDEK) with prescribed integral curvature $\cK=\kappa$ 
given in (\iCURV) and finite area $\cA=\alpha$ given in (\AREA).
	Looking only for radial solutions reduces our 
PDEs for $K$ and $u$ to two ODEs.
	We transform these ODEs for $K$ and $u$ into a potential
scattering problem for a single Newtonian particle in $\RR^2$ and
solve this scattering problem by fixed point arguments aided with
gradient flow techniques.
	This strategy is adapted from [\KieJllPHPL] where a 
different self-consistent Gauss curvature problem is considered.

	Let $\xi=f_\xi(t),\ \eta=f_\eta(t)$ be the time-dependent 
Cartesian coordinates of a point in $\RR^2$ which moves according 
to the Newtonian equations of motion
$$
\eqalignno{
&
	{\dd^2\xi\over \dt^2} 
= 
	- {\partial V\over \partial \xi} \, ,
&\eqlbl\NewtonXI\cr
&
	{\dd^2\eta\over \dt^2}  
= 
	- {\partial V\over \partial \eta} \, ,
&\eqlbl\NewtonETA\cr}
$$
in a fixed external potential 
$$
V(\xi,\eta) = {1\over 2} \eta e^{2\xi}.
\eqno\eqlbl\NewtonPOT
$$
	We will sometimes write $\xi(t)$, $\eta(t)$ and $\dot\xi(t)$,
$\dot\eta(t)$ to  denote solutions and their time derivatives. 
	We seek solutions of (\NewtonXI), (\NewtonETA), (\NewtonPOT) 
that satisfy the asymptotic conditions 
$$
\eqalignno{
\lim_{t\to -\infty} {\xi(t) - t}	&= \xi_{\rm in}
&\eqlbl\ACone\cr
\lim_{t\to -\infty} {\eta(t) } 		&= \eta_{\rm in}
&\eqlbl\ACtwo}
$$
for suitable real constants $\xi_{\rm in}$ and $\eta_{\rm in}$ such
that there exists a $\Theta\in (-\pi,-\pi/2)$ such that
$$
\eqalignno{
\lim_{t\to + \infty} {\xi(t) \over t} 	&= \cos\Theta,
&\eqlbl\ACthree\cr
\lim_{t\to + \infty} {\eta(t) \over t} 	&= \sin\Theta.
&\eqlbl\ACfour}
$$
	Clearly, the asymptotic conditions 
(\ACone), (\ACtwo), (\ACthree), (\ACfour) imply that 
asymptotically in the infinite past and the infinite future
the particle performs a linear, unaccelerated  motion.
	These two ``asymptotically free motions'' are linked
by a deflection of the particle off of its initial direction by 
an angle $\Theta$, which is effected by the external potential $V$. 
	Our problem thus belongs in the category ``potential scattering.''
\smallskip

\noindent
{\bf Theorem 3.1:} {\it 
	For each $\Theta\in (-\pi,-\pi/2)$
there exists a constant $\eta_{\rm in} >0$,  
such that for each $\xi_{\rm in}\in \RR$ there exists a unique 
solution pair $\xi(t),\ \eta(t)$ of 
(\NewtonXI), (\NewtonETA), (\NewtonPOT) 
satisfying (\ACone), (\ACtwo), (\ACthree), (\ACfour). 
	Within the family of solutions belonging to the same 
$\eta_{\rm in}$ we can switch from one solution 
to another by means of the transformation  
$\xi_{\rm in}\to\xi_{\rm in}^\pr$ combined with a corresponding  
time translation $t\to t + \xi_{\rm in} -\xi_{\rm in}^\pr$.
	This transformation leaves $\Theta$ unchanged.}
\smallskip

	Before we prove Theorem 3.1, we first show that our
Theorem 1.2 is a corollary of Theorem 3.1.
\smallskip

{\it Proof of Theorem 1.2}: 
	Let $\xi=f_\xi(t),\eta = f_\eta(t)$ denote the motion of a 
Newtonian point particle in $\RR^2$ according to (\NewtonXI), (\NewtonETA)
with  $V$ given in (\NewtonPOT), having asymptotic behavior given by 
(\ACone), (\ACtwo), (\ACthree), (\ACfour).
	By Theorem 3.1, such a motion exists. 
	Inserting (\NewtonPOT) into (\NewtonXI) and (\NewtonETA), 
the equations of motion read explicitly
$$
\eqalignno{
&
	{\dd^2\xi\over \dt^2} 
= 
	- \eta e^{2\xi}\, ,
&\eqlbl\NewtonXInew\cr
&
	{\dd^2\eta\over \dt^2}  
= 
	- {1\over 2} e^{2\xi}\, .
&\eqlbl\NewtonETAnew\cr}
$$
	We now set $t = \ln r$ for $r>0$, define
$$
\ol{u}(r) =  f_\xi(\ln r) - \ln r -{1\over 4}\ln 2
\eqno\eqlbl\Udefine
$$
and
$$
\ol{K}(r) = \sqrt{2}\, f_\eta(\ln r)\, ,
\eqno\eqlbl\Kdefine
$$
and find that for $r>0$, the functions $\ol{u}(r)$ and $\ol{K}(r)$ satisfy
$$
-{1\over r }{{\rm d} \over \dr} \left( r {{\rm d} \over \dr} \ol{u}(r)\right)
= \ol{K}(r) e^{2\ol{u}(r)}
\eqno\eqlbl\Newtonu
$$
and 
$$
-{1\over r }{{\rm d} \over \dr} \left( r {{\rm d} \over \dr} \ol{K}(r)\right)
= e^{2\ol{u}(r)}.
\eqno\eqlbl\NewtonK
$$
	Moreover, we can set $\ol{K}(0)=\sqrt{2}\,\eta_{\rm in}$, and 
$\ol{u}(0) = \xi_{\rm in} -{1\over 4}\ln 2 $. 
	Identifying $r = |x-x_*|$ for $x\in \RR^2$, with
$x_*$ the arbitrary center of symmetry,  we recognize that
(\Newtonu) is (\PDEu), and (\NewtonK) is (\PDEK), for  radially 
symmetric $K(x)=\ol{K}(|x-x_*|)$ and $u(x)= \ol{u}(|x-x_*|)$.
	Furthermore, from (\NewtonETAnew) it follows that $K(x)$ is
decreasing away from $x_*$, and from (\ACtwo) we have 
$K^*=\sqrt{2}\, \eta_{\rm in}$. 
	From (\ACfour) it follows that 
$K(x) \sim -\sqrt{2}\,\sin\Theta\,\ln |x|$
as $|x|\to\infty$, as claimed. 
	We have the identification $2\pi\sqrt{2}\, \sin\Theta = \alpha$, 
so that from (\ACfour) and (\Newtonu) it follows that $\cK(u)
=\left( 4\pi+\sqrt{16\pi^2 -2\alpha^2}\,\;\right)\!/2 \in (2\pi,4\pi)$, 
as demanded by (\pohoALT).
	Finally, translations $t\mapsto t +t_0$ combined with an 
associated translation $\xi\to \xi + \xi_0$ correspond to scalings
$r\mapsto kr$ combined with translations $u \mapsto u - \ln k$, which 
together with the indeterminacy of $x_*$ proves that $u$ is unique modulo 
the conformal transformations listed in Theorem 1.2.\qed
\smallskip

	It remains to prove Theorem 3.1.
\vfill\eject

	We begin by listing the symmetries of the ODE system
(\NewtonXInew), (\NewtonETAnew), which are:

\smallskip
{\it $\bullet$ Invariance under time translations $t\to t + t_0$;}

\smallskip
{\it $\bullet$ Invariance under time reversal $t\to -t$;}

\smallskip
{\it $\bullet$ Invariance under the homologous
transformations $\xi\to\xi +\xi_h$ and $t\to e^{-\xi_h}t$.}
\smallskip

\noindent
	By E. Noether's theorem, invariance under time translations
is associated with the conservation law for the total (kinetic plus 
potential) energy $E$ of the Newtonian unit mass point, where
$$
2E = \dot\xi^2 + \dot\eta^2 + \eta e^{2\xi} .
\eqno\eqlbl\COM
$$ 
	Under the homologous transformations $\xi\to\xi +\xi_h$ 
and $t\to e^{-\xi_h}t$ the conserved quantity $E$ transforms as 
$E\to e^{2\xi_h}E$. 
	Hence, to obtain all solutions of (\NewtonXInew), 
(\NewtonETAnew) it suffices to obtain all solutions for three
generic values of $E$, say $E=E_+>0$, $E=0$, and $E=E_- <0$. 
	For the motion of interest to us, the asymptotic conditions
(\ACone) and (\ACtwo) give
$$
E=1/2 \, .
\eqno\eqlbl\COMvalue
$$

\noindent
{\bf Lemma 3.2:} {\it A solution $\xi=f_\xi(t),\ \eta = f_\eta(t)$ of the
equations of motion (\NewtonXInew), (\NewtonETAnew)
satisfying (\ACone)--(\ACfour) is restricted to
the region $\{(\xi,\eta)\in \RR^2:\, \eta < e^{-2\xi}\}$.}
\smallskip

{\it Proof of Lemma 3.2}: 
	Clearly, since the kinetic energy is non-negative,
(\COMvalue) cannot be achieved in the ``$E=1/2$ forbidden zone'' 
where $\eta > e^{-2\xi}$.
	Hence, a solution $\xi=f_\xi(t),\ \eta = f_\eta(t)$ of the
Newtonian equations of motion (\NewtonXInew), (\NewtonETAnew)
satisfying (\ACone)--(\ACfour) is confined to
the region $\{(\xi,\eta)\in \RR^2:\, \eta \leq e^{-2\xi}\}$.
	It remains to show that a solution cannot have a point
in common with the boundary $\{\eta = e^{-2\xi}\}$ of the $E=1/2$ 
forbidden zone.

	The boundary $\eta = e^{-2\xi}$ of the $E=1/2$ forbidden zone 
consists of all points $(\xi,\eta)\in \RR^2$ for which 
$E=1/2$ is achieved iff $\dot\xi = 0 = \dot\eta$. 
	Recall that a {\it singular point} on a trajectory is a point
at which both $\dot\xi = 0$ and $\dot\eta =0$; hence, the boundary of 
the  $E=1/2$ forbidden zone consists of all the possible singular
points.
	A trajectory which contains (at least one) singular point
is called a {\it singular trajectory}. 
	Thus, a singular trajectory has at least one point in
common with the boundary of the $E=1/2$ forbidden zone.
	On the other hand, it follows immediately from (\NewtonETAnew) 
that there can be at most one singular point on a singular trajectory,
hence a singular trajectory has exactly one singular point.
	By the time translation invariance of (\NewtonXInew), (\NewtonETAnew) 
we can assume that this point is reached at $t=0$. 
	By the time reversal  invariance of (\NewtonXInew), (\NewtonETAnew) 
it now follows that on a singular trajectory the forward motion with respect
to $t=0$ is identical to the backward motion with respect to $t=0$.
	This in turn implies that the asymptotic conditions are
symmetric under time-reversal as well.
	But then by (\ACone) and (\ACfour) we conclude that
$\cos\Theta = 1$, which implies $\sin\Theta = 0$,
which contradicts the condition that $\Theta\in (-\pi,-\pi/2)$. 
	Hence, the motion on a singular trajectory cannot satisfy all 
our asymptotic conditions.
	Put differently, a solution to our equations of motion which 
does satisfy all asymptotic conditions cannot be singular. 
	Our Lemma 3.2 is proved.\qed
\vfill\eject

\noindent
{\bf Lemma 3.3:}
{\it 	
	Let $\xi=f_\xi(t),\ \eta =f_\eta(t)$ 
solve (\NewtonXInew), (\NewtonETAnew)
for the asymptotic conditions (\ACone), (\ACtwo).
	Then the map $f = f_\xi\circ f_\eta^{-1}$ is well defined on the set 
$f_\eta(\RR)$, and we have $\xi= f(\eta)$. 
	Furthermore, there exists a unique $\eta_\sim < \eta_{\rm in}$ such 
that $f$ is strictly convex for $\eta < \eta_\sim$ and strictly concave for 
$\eta > \eta_\sim$.
}
\smallskip

{\it Proof of Lemma 3.3}: 
	By integrating (\NewtonETAnew) once, using 
(\ACtwo), we have
$$
	\dot\eta(t) 
=  
	-{1\over 2} \int_{-\infty}^t e^{2\xi(s)}\dd s .
\eqno\eqlbl\etadotINT
$$	
	Clearly, the map $t\mapsto \dot\eta(t)$ is strictly negative 
for all $t>-\infty$; hence, the map $t\mapsto \eta= f_\eta(t)$ is strictly
monotonically decreasing and thus invertible, giving $t=f_\eta^{-1}(\eta)$.

	Next, let $^\pr$ denote derivative with respect to $\eta$.  
	Along a trajectory $\xi=f_\xi(t),\ \eta =f_\eta(t)$ that solves
(\NewtonXInew), (\NewtonETAnew) for the asymptotic conditions (\ACone), 
(\ACtwo), we then have  
$$
	f^\ppr(\eta) = {{\rm d}^2 \xi\over {\rm d}\eta^2} 
= 
	{1\over \dot\eta{}^3}
      \left(\ddot \xi\dot \eta -\ddot \eta\dot \xi\right) 
= 
	{\exp (2 \xi) \over 2 \dot \eta{}^3} 
	\left(\dot\xi  - 2\eta\dot \eta \right),
\eqno\eqlbl\trajectddot
$$
the middle and right sides evaluated at $t$, the 
left side at $\eta = f_\eta(t)$. 
	By (\etadotINT), the map $t\mapsto \dot\eta^3$ is negative and
strictly monotonically decreasing.
	Next notice that by multiplying (\NewtonETAnew) by  $2\eta$ and 
subtracting that result from (\NewtonXInew) we get 
$$
{{\rm d}^2 \xi\over {\rm d}t^2} - 2\eta{{\rm d}^2\eta\over {\rm d}t^2} =0\, .
\eqno\eqlbl\ddotsubtraction
$$
	Upon integrating (\ddotsubtraction) from $-\infty$ to $t$, using
integration by parts, we obtain
$$
\left( \dot\xi  - 2\eta\dot \eta \right)  (t)
= 1 -2 \int_{-\infty}^t \dot{\eta}^2 (s) \ds \, .
\eqno\eqlbl\XYdot
$$
	Since $t\mapsto \dot\eta^2(t)$ is positive and 
strictly monotonically increasing, by (\XYdot) we now conclude
that the map $t\mapsto \int_{-\infty}^t \dot{\eta}^2 (s) \ds$ is  
strictly monotonically increasing and strictly convex.
	Therefore there exists a unique $t_\sim$ such that the \rhs\ of 
(\XYdot) is strictly positive for $t<t_\sim$ and strictly negative 
for $t>t_\sim$.
	Setting $\eta_\sim \equiv f_\eta(t_\sim)$, we then conclude 
that the \rhs\ of (\XYdot) evaluated at $t=f_\eta^{-1}(\eta)$ is strictly 
positive for $\eta >\eta_\sim$ and strictly negative for $\eta < \eta_\sim$.
	We thus conclude from (\trajectddot) that along the trajectory 
$\xi = f(\eta)$ we have 
$$
f^\ppr(\eta) \ \cases{ \ > 0 \ {\rm for}\ \eta <\eta_\sim\cr\cr
			\ < 0 \ {\rm for}\ \eta >\eta_\sim \ \ ,}
\eqno\eqlbl\convexityTRAJECTORY
$$
as claimed.\qed
\smallskip

	By the convexity of $\eta\mapsto \xi = f(\eta)$ for 
$\eta<\eta_\sim$ 
it follows that a solution $\xi=f_\xi(t),\ \eta =f_\eta(t)$ of 
(\NewtonXInew), (\NewtonETAnew), (\ACone), (\ACtwo) which 
satisfies a linear bound $f(\eta) < A\eta +B$ for some constants 
$A >0$ and $B$ necessarily satisfies the asymptotic conditions 
(\ACthree), (\ACfour) for some $\Theta\in (-\pi,-\pi/2)$.
	Part of our existence proof will concentrate on proving
that for $\eta_{\rm in}$ large enough such a linear bound on $f$ 
exists.
\vfill\eject

	On the other hand, such a linear bound on $f$ will fail to
exist if $\eta_{\rm in}$ is negative. 
	Namely, by (\etadotINT) we have $\dot\eta(t) < 0$ 
for all $t>-\infty$, which implies that 
$\sup_t \eta(t) = \lim_{t\to-\infty}\eta(t)$.
	By (\ACtwo) we then have $\sup_t\eta(t) = \eta_{\rm in}$.
	Therefore, if $\eta_{\rm in} \leq 0$, we conclude that
$\eta(t)<0$ for all $t>-\infty$. 
	Integrating now (\NewtonXInew) once, using (\ACone),  we obtain 
$$
	\dot\xi(t) 
= 
	1 - \int_{-\infty}^t \eta(s)e^{2\xi(s)}\dd s.
\eqno\eqlbl\xidotINT
$$
	Since $\eta(t) <0$ for all $t>-\infty$ if $\eta_{\rm in}\leq 0$, 
(\xidotINT) now implies that $\dot\xi(t) > 0$ for all $t>-\infty$, 
which contradicts the asymptotic condition (\ACfour), which is negative 
for $\Theta\in (-\pi,-\pi/2)$. 
	Hence, we have proven 
\smallskip
\noindent
{\bf Proposition 3.4:} 
{\it
	If a solution $\xi=f_\xi(t),\ \eta = f_\eta(t)$ of 
(\NewtonXInew), (\NewtonETAnew) satisfies (\ACone)---(\ACfour), with 
$\Theta\in (-\pi,-\pi/2)$, then $\eta_{\rm in} >0$.
}
\smallskip

	Next, let $T=T(\xi_{\rm in},\eta_{\rm in})$ be
the instant where the maximal Cauchy development terminates. 
	Then for $t<T$ the system of differential equations
(\NewtonXInew), (\NewtonETAnew) with asymptotic conditions (\ACone), (\ACtwo)
is equivalent to the coupled system of nonlinear integral equations
$$
\eqalignno{
	\xi(t)
& = 
	\xi_{\rm in}+t  - \int_{-\infty}^t\int_{-\infty}^s 
		\eta(\ts)e^{2\xi(\ts)}\dts\ds ,
&\eqlbl\xiINT\cr
	\eta(t) 
&= 
	\eta_{\rm in} -{1\over 2} \int_{-\infty}^t\int_{-\infty}^s
	e^{2\xi(\ts)}\dd \ts\dd s,
&\eqlbl\etaINT}
$$
obtained by integrating (\xidotINT) using (\ACone), and integrating
(\etadotINT), using (\ACtwo).
	We remark that there do exist solutions that blow up at a
finite time $T<\infty$ if $\eta_{\rm in}$ is below some critical 
value (in particular, this is the case if $\eta_{\rm in} <0$). 

	To analyze (\xiINT), (\etaINT), we study
the coupled iteration sequences 
$$
\eqalignno{
\hskip-1truecm 
	\xi^{(n)}(t)
& = 
	\xi_{\rm in} +t  - \int_{-\infty}^t\int_{-\infty}^s 
		\eta^{(n)}(\ts)e^{2\xi^{(n)}(\ts)}\dts\ds ,
&\eqlbl\iteraXI\cr
	\eta^{(n+1)}(t)
& = 	
	\eta_{\rm in} -{1\over 2} \int_{-\infty}^t\int_{-\infty}^s 
		\,e^{2\xi^{(n)}(\ts)}\dts\ds\, ,
&\eqlbl\iteraETA}
$$
$n\geq 0$, with the starting function $\eta^{(0)}$ given by
$$
	\eta^{(0)}(t)
 \equiv
	\eta_{\rm in}.
\eqno\eqlbl\iteraETAnull
$$
	By inspection one readily checks that, if the iteration sequences
(\iteraXI), (\iteraETA) with starting function
(\iteraETAnull) converge for all $t<T$, then they converge to
functions $\xi=f_\xi(t),\ \eta = f_\eta(t)$ solving
(\xiINT), (\etaINT).
	We have to show that for large enough $\eta_{\rm in}$, 
the sequences converge to functions satisfying also
(\ACthree) and (\ACfour), in which case $T=\infty$. 
\smallskip

\noindent
{\bf Lemma 3.5:} 
{\it 
	For $\eta_{\rm in} >0$, the maps $n\mapsto \xi^{(n)}$ and
$n\mapsto \eta^{(n)}$ defined jointly by the iteration sequences 
(\iteraXI), (\iteraETA) with starting function (\iteraETAnull) are 
pointwise increasing, respectively decreasing, for each fixed 
$t> -\infty$. 
}
\smallskip

{\it Proof of Lemma 3.5}: 
	The claim of Lemma 3.5 follows by standard 
sub- and supersolution techniques.
	Using (\iteraETAnull) we see that (\iteraXI) for $n=0$ reads
$$
	\xi^{(0)}(t) 
 = 
	\xi_{\rm in} +t - \eta_{\rm in} \int_{-\infty}^t\int_{-\infty}^s 
		 e^{2\xi^{(0)}(\ts)}\dts\ds .
\eqno\eqlbl\iteraXInull
$$
	For $\eta_{\rm in} >0$ the nonlinear integral equation 
(\iteraXInull) is solved uniquely by
$$
	\xi^{(0)}(t) 
=  
	-\ln\cosh \big(t + \xi_{\rm in} - \ln(2/\sqrt{\eta_{\rm in}})\big)
	- \ln \sqrt{\eta_{\rm in}}\, .
\eqno\eqlbl\XInull
$$
	Thus, for all $p>0$ and for all $t>-\infty$ the integral 
$\int_{-\infty}^t\int_{-\infty}^s |\ts|^p e^{2\xi^{(0)}(\ts)}\dts\ds$
exists; in particular, the integral exists for $p=0$. 
	Therefore (\iteraETA) for $n=0$ is well defined for all 
$t>-\infty$, and by integration we find $\eta^{(1)}(t)$ to be given by
$$ 
	\eta^{(1)}(t) 
=  
	- {1\over 2\eta_{\rm in}} 
	  \ln\cosh \big(t +\xi_{\rm in} -\ln(2/\sqrt{\eta_{\rm in}})\big)
	- {t\over 2\eta_{\rm in}}
	- {\xi_{\rm in}\over 2\eta_{\rm in}} 
	+ \eta_{\rm in} 
	-  {1\over 4\eta_{\rm in}}\ln \eta_{\rm in} .
\eqno\eqlbl\ETAone
$$
	Clearly, $\eta^{(1)}(t)\to \eta_{\rm in}$ as $t\to -\infty$, 
and $\eta^{(1)}(t)\sim - {1\over\eta_{\rm in}} t$ as $t\to +\infty$;
moreover, $\eta^{(1)}(t) < \eta_{\rm in} = \eta^{(0)}$ for all $t$, 
which is seen by inspection of (\ETAone) but also 
follows immediately from (\iteraETA).
	Hence, (\iteraXI) with $n=1$ has a well defined 
solution $\xi^{(1)}(t)$ for all $t < T^{(1)}$. 
	Moreover, (\iteraXI) implies at once that 
$\xi^{(1)}(t) > \xi^{(0)}(t)$ for all $t$ for which $\xi^{(1)}$ exists. 
	Hence, we conclude that $\eta^{(2)}<\eta^{(1)}$, and so on
by induction.\qed
\smallskip

\noindent
{\bf Lemma 3.6:} 
{\it 
Let $\xi^{(n)}(t),\ \eta^{(n)}(t)$ solve (\iteraXI) (\iteraETA),
(\iteraETAnull). 
	Then there exists a $T_0=T_0(\xi_{\rm in},\eta_{\rm in})$,
independent of $n$, satisfying the bound 
$$
T_0 >  \ln(2\sqrt{2\eta_{\rm in}}) - \xi_{\rm in},
\eqno\eqlbl\TNULLlowerBOUND
$$ 
such that for all $t< T_0$ and for all $n$ we have 
$$
\eta^{(n)}(t) > 0
\eqno\eqlbl\ETAlowerBOUND
$$
and
$$
\xi^{(n)}(t) < \xi_{\rm in} +t.
\eqno\eqlbl\XIupperBOUND
$$
} 
\smallskip
{\it Proof of Lemma 3.6}:  
	Clearly, for each $n$ the function  $t\mapsto \eta^{(n)}(t)$
is strictly monotonic decreasing and strictly concave.
	Since $\eta_{\rm in}>0$, there exists a 
unique $T^{(n)}_0(\xi_{\rm in},\eta_{\rm in})$ such that
$\eta^{(n)}\big(T^{(n)}_0\big)=0$.
	Moreover, since the iteration map $n\mapsto \eta^{(n)}(t)$ is 
decreasing for each $t$, we conclude that the sequence 
$n\mapsto T^{(n)}_0(\xi_{\rm in},\eta_{\rm in})$ is decreasing, too.
	We need to show that it has a lower bound $T_0 > -\infty$.
 
	Now, by what we just said, it follows with (\iteraXI) that
for all $t < T^{(n)}_0$ we have the $n$-independent upper bound 
(\XIupperBOUND) for $\xi^{(n)}(t)$.
	This in turn implies that for all $t < T^{(n)}_0$
we have the $n$-independent lower bound
$$
\eta^{(n)}(t) > \eta_{\rm in} - {1\over 8}e^{2\xi_{\rm in} +2t}.
\eqno\eqlbl\ETAlowerBOUNDt
$$
	By setting the r.h.s. of (\ETAlowerBOUNDt) equal to zero
we obtain the $n$-independent lower bound r.h.s.(\TNULLlowerBOUND)
valid for all $T_0^{(n)}$; thus the $T_0^{(n)}$ are bounded below 
independently of $n$ by some $T_0$  satisfying (\TNULLlowerBOUND),
and our Lemma follows at once.\qed
\smallskip

\noindent
{\bf Corollary 3.7:} {\it 
The sequence $n\mapsto \big(\xi^{(n)}(t),\eta^{(n)}(t)\big)$ defined
by (\iteraXI) (\iteraETA), (\iteraETAnull) converges pointwise for
all $t<T$ (the  life span of the maximal Cauchy development)
to a solution $\big(\xi_*(t),\eta_*(t)\big)$ of (\xiINT) and (\etaINT),
and this is the unique solution to (\NewtonXInew), (\NewtonETAnew), 
satisfying (\ACone) and (\ACtwo).}
\smallskip

{\it Proof of Corollary 3.7}:  
	By Lemma 3.5, the sequence
$n\mapsto \big(\xi^{(n)}(t),\eta^{(n)}(t)\big)$ defined
by (\iteraXI) (\iteraETA), (\iteraETAnull) is pointwise increasing
for $\xi$ and decreasing for $\eta$.
	By Lemma 3.6, for all $t<T_0$ the $\xi$ sequence is bounded
above and the $\eta$ sequence bounded below independently of $n$.
	Hence,  these two sequences converge for $t<T_0$ to solutions
$\xi_*(t)$ and $\eta_*(t)$ of (\xiINT) and (\etaINT).
	Furthermore, by our sharp upper and lower bounds on any
solution $\xi(t)$ and $\eta(t)$ for  $t<\tau\ll T_0$, we can 
easily show that the fixed point map defined by  (\xiINT) and 
(\etaINT) is a contraction mapping in the set of integrable 
functions on $(-\infty,\tau)$ equipped with exponentially 
weighted $L^1$ norm, hence the solutions
$\xi_*(t)$ and $\eta_*(t)$ of (\xiINT) and (\etaINT) are
unique for $t<\tau$. (We skip the details of the contraction mapping
proof here because below we reprove the uniqueness by a different
argument that will be needed in the sequel.)

	Next, we can now pick any particular $t_{0} < \tau$ as 
new initial time and solve (\NewtonXInew), (\NewtonETAnew) for
$t> t_{0}$ as regular initial value problem with data 
$\xi_*(t_{0})$ and $\eta_*(t_{0})$.
	Standard ODE results now guarantee that this initial 
value problem has a unique solution for all $t\in (t_0,T)$, 
and this solution satisfies (\xiINT) and (\etaINT) and moreover
can be computed with (\iteraXI), (\iteraETA), (\iteraETAnull).
	Thus, the solution $\big(\xi_*(t),\eta_*(t)\big)$
is continued uniquely  from $t\in(-\infty, t_0]$ to $t\in (t_0,T)$,
and this proves the corollary.\qed
\smallskip

	Having a unique solution to  (\xiINT) and (\etaINT) for all 
$t<T$, where by uniqueness we now also know that 
$T=T(\xi_{\rm in},\eta_{\rm in})$, we can bootstrap to a sharper 
upper bound on $\xi(t)$.
\smallskip

\noindent
{\bf Lemma 3.8:} 
{\it Let  $\big(\xi(t),\eta(t)\big)$  solve (\NewtonXInew), (\NewtonETAnew) 
for the asymptotic conditions (\ACone), (\ACtwo). 
	Let  $T_{1/2}$ be defined by $\eta(T_{1/2}) = \eta_{\rm in}/2$.
	Then, for $T_{1/2}$ we have the lower bound
$$
T_{1/2} >\ln( 2\sqrt{\eta_{\rm in}}) - \xi_{\rm in},
\eqno\eqlbl\THALFlowerBOUND
$$
and for all $t\in (-\infty,T_{1/2})$ we have the upper bound 
$\xi(t) < \hat\xi(t)$, where 
$$
	\hat\xi(t) 
=
	-\ln\cosh \big(t + \xi_{\rm in} - \ln(2\sqrt{2/\eta_{\rm in}})\big)
	- \ln \sqrt{\eta_{\rm in}/2}\, .
\eqno\eqlbl\XIhat
$$}

{\it Proof of Lemma 3.8:} 
	As for $T_{1/2}$, for all $t<T_{1/2}$ we have the
lower bound (\ETAlowerBOUNDt) for $\eta$.
	By setting the r.h.s. of (\ETAlowerBOUNDt) equal to $\eta_{\rm in}/2$,
we obtain the lower bound (\THALFlowerBOUND).

	Since $\eta(t) > \eta_{\rm in}/2$ for
$t< T_{1/2}$, we find from (\xiINT) that the solution to 
$$
	\hat\xi(t)
=
	\xi_{\rm in}+t  - {1\over 2}\eta_{\rm in}
		\int_{-\infty}^t\int_{-\infty}^s e^{2\hat\xi(\ts)}\dts\ds 
\eqno\eqlbl\xiINTsuperEQ
$$
is a supersolution  for $\xi(t)$  for all $t< T_{1/2}$.
	For $\eta_{\rm in} >0$ the nonlinear integral equation 
(\xiINTsuperEQ) is solved uniquely by (\XIhat).\qed
\smallskip

\noindent
{\bf Lemma 3.9:} 
{\it There exists some $\eta_{\rm in}^{\rm crit} >0$ 
such that when $\eta_{\rm in} >\eta_{\rm in}^{\rm crit}$, 
then $\xi(t)$ has a maximum at some finite $T_M < T_0$
(the same $T_0$ as in Lemma 3.6).
	In that case, at $t=T_0$ we have  the bounds
$$
	\xi(T_0) 
< 	
	-\ln\cosh \ln\big(\eta_{\rm in}/\sqrt{2}\big)
	- \ln \sqrt{\eta_{\rm in}/2}\, ,
\eqno\eqlbl\xiNULLbound
$$
$$
\dot\xi(T_0) 
< 
{-\ln\cosh \ln\big(\eta_{\rm in}/\sqrt{2}\big)+ \ln \sqrt{2}
\over
\ln{\eta_{\rm in}} } <0,
\eqno\eqlbl\xiDOTbound
$$
and 
$$
\dot\eta(T_0) 
> - \sqrt{1 - \left({\ln\cosh \ln\big(\eta_{\rm in}/\sqrt{2}\big)-\ln\sqrt{2}
			\over
			\ln\eta_{\rm in}}  \right)^2}\, .
\eqno\eqlbl\etaDOTbound
$$
}

{\it Proof of Lemma 3.9}: 
	The proof exploits the convexity properties of $\xi(t)$ for $t> T_0$.
	Namely, by (\NewtonXInew), for all $t>T_0$, $\xi(t)$ is concave (i.e. 
convex down).
	Furthermore, for all $t\in (-\infty, T_{1/2})$ (recall that
$T_{1/2}< T_0$), $\xi(t)$ satisfies the manifestly concave sandwich bounds 
$\xi^{(0)}(t) <\xi(t) < \hat\xi(t)$, given by (\XInull) and (\XIhat).
	Next, let $T_M^{(0)}$  and $\hat{T}_M$ be the instants at which 
$\xi^{(0)}(t)$ and $\hat\xi(t)$ take their respective maximum, and
let $\overline{T}_{1/2}$ be given by the r.h.s. of (\THALFlowerBOUND).
	It is readily seen that 
$T_M^{(0)} =  \ln(2/\sqrt{\eta_{\rm in}}) - \xi_{\rm in}$ and
$\hat{T}_M =  \ln(2\sqrt{2/\eta_{\rm in}}) - \xi_{\rm in}$. 
	For $\eta_{\rm in}> \sqrt{2}$ we have the ordering
$-\infty < T_M^{(0)} < \hat{T}_M < \overline{T}_{1/2} < {T}_{1/2} < T_0$.
	Furthermore, we have the monotonic behavior that, 
as $\eta_{\rm in}\nearrow$, we have $T_M^{(0)}\searrow$ and 
$\hat{T}_M\searrow$, but $\overline{T}_{1/2}\nearrow$.
	Now let $\wtilde{\eta}_{\rm in}^{\rm crit}$ be the unique 
solution of $\xi^{(0)}(T_M^{(0)}) = \hat\xi(\overline{T}_{1/2})$.
	After a simple manipulation, we see that 
$\wtilde{\eta}_{\rm in}^{\rm crit}$ is given by
$$
{\wtilde{\eta}_{\rm in}^{\rm crit}}
=
	\sqrt{2} \exp\arcosh 2 .
\eqno\eqlbl\etaCRITbound
$$
	Clearly, ${\wtilde{\eta}_{\rm in}^{\rm crit}}>\sqrt{2}$.
	But then, by the geometry of the concave sandwich bounds
and the ordering and monotonic behavior of the various instances of time, 
we conclude that for all $\eta_{\rm in} > \wtilde{\eta}_{\rm in}^{\rm crit}$
we have that $\xi(T_M^{(0)}) > \xi(\overline{T}_{1/2})$, and 
therefore $\xi(t)$ has a unique maximum at some 
$T_M < \overline{T}_{1/2})$ whenever 
$\eta_{\rm in} > \wtilde{\eta}_{\rm in}^{\rm crit}$.

	Next, whenever $\eta_{\rm in} > \wtilde{\eta}_{\rm in}^{\rm crit}$
so that $\xi(t)$ has a maximum for $T_M < T_0$, it follows directly
from (\NewtonXInew) that $\dot\xi(t) < 0$ for all $T_M < t < T_0$.
	Therefore, we conclude that $\xi(T_0) < \hat\xi(\overline{T}_{1/2})$,
and this gives the bound (\xiNULLbound).

	The bound (\xiDOTbound) follows once again by convexity arguments.
	Namely, by the concavity of $\xi(t)$ for $t>T_0$, it follows
that whenever $\eta_{\rm in} > \wtilde{\eta}_{\rm in}^{\rm crit}$,
we have that $\dot\xi(T_0) < \dot\xi(\overline{T}_{1/2})$.
	To estimate $\dot\xi(\overline{T}_{1/2})$ we simply compute 
the slope of the straight line joining 
the maximum of $\xi^{(0)}$ with $\hat\xi(\overline{T}_{1/2})$.
	By the convexity of these sandwich bounds on $\xi$ 
it follows right away that the slope of that straight line dominates
$\dot\xi(\overline{T}_{1/2})$. 
	This is the content of (\xiDOTbound).

	Finally, at $t= T_0$ we have $\eta(T_0)=0$, so that by
the energy law (\COMvalue)  we have that 
$\dot\xi(T_0)^2 + \dot\eta(T_0)^2 =1$. 
	But $\dot\eta(t) < 0$ for all $t$, hence at $t=T_0$ we
have $\dot\eta(T_0) = - \big(1 - \dot\xi(T_0)^2\big)^{1/2}$.
	With (\xiDOTbound) we now obtain (\etaDOTbound).
	Finally, from the way it is constructed it is manifestly clear 
that $\wtilde{\eta}_{\rm in}^{\rm crit}$ is an upper
estimate for ${\eta}_{\rm in}^{\rm crit}$.\qed

	We now turn to the time zone $t\geq T_0$ 
and derive an asymptotically linear upper bound for 
$\xi(t)$ and an asymptotically linear lower bound for $\eta(t)$,
valid whenever $\eta_{\rm in} > \wtilde{\eta}_{\rm in}^{\rm crit}$.
	Thus, $\eta_{\rm in} > \wtilde{\eta}_{\rm in}^{\rm crit}$, 
and let $\eps \ll 1$. 
	For $t\geq T_0$ define two maps
$F_\eps$ and $G_\eps$ from $C^0\times C^0$ to $C^0$ by
$$
\eqalignno{
\hskip-1truecm 
	F_\eps(X,Y)(t) 
& = 	X(t) - \eps
	\left(
		X(t) -\dot{X}(T_0) (t-T_0) - X(T_0) 
		+ \int_{T_0}^t\int_{T_0}^s Y(s^\pr)e^{2X(s^\pr)}\dss\ds 
	\right),
&\cr
& &\eqlbl\iteraF\cr
	G_\eps(X,Y)(t) 
& = 	Y(t) -\eps
	\left( 
		Y(t) -\, \dot{Y}(T_0) (t-T_0)
	    + \int_{T_0}^t\,\int_{T_0}^s\,{1\over 2}\,e^{2X(s^\pr)}\dss\ds\,
	 \right),
&\eqlbl\iteraG}
$$
where $t\mapsto X(t)$ and $t\mapsto Y(t)$ are any two continuous
functions that satisfy the initial bounds 
$X(T_0)< $ r.h.s.(\xiNULLbound), 
$\dot{X}(T_0) <$r.h.s.(\xiDOTbound),
$Y(T_0)=0$, and 
r.h.s.(\etaDOTbound)$<\dot{Y}(T_0)<0$. 
	Now consider the coupled iteration sequences
$$
\eqalignno{
\hskip-1truecm 
	X^{(n+1)}(t)
& = 
	F_\eps\big(X^{(n)},Y^{(n)}\big),
&\eqlbl\iteraX\cr
	Y^{(n+1)}(t)
& = 	
	G_\eps\big(X^{(n)},Y^{(n)}\big),
&\eqlbl\iteraY}
$$
with the starting functions
$$
\eqalignno{
	X^{(0)}(t) 
& = 
	\dot{X}(T_0) (t-T_0)+ X(T_0); \qquad t\geq T_0
&\eqlbl\iteraXnull\cr
	Y^{(0)}(t)
& = 
	\dot{Y}(T_0) (t-T_0);\qquad\qquad\qquad t\geq T_0.
&\eqlbl\iteraYnull}
$$
\smallskip

\noindent
{\bf Lemma 3.10:} {\it 
	The maps $n\mapsto X^{(n)}$ and $n\mapsto Y^{(n)}$ defined
jointly by the iteration sequences (\iteraX), (\iteraY) 
with (\iteraF), (\iteraG)  and starting functions 
(\iteraXnull), (\iteraYnull) are increasing, respectively 
decreasing, pointwise for all $t>T_0$. 
}
\smallskip

{\it Proof of Lemma 3.10}:  
	We  prove Lemma 3.10 by induction. 
	
	First, we obviously have $Y^{(1)}(t)<Y^{(0)}(t)$ for all $t>T_0$.  
	Since $\dot{Y}(T_0)<0$ by (\etadotINT), 
we also have $Y^{(0)}(t)<0$
for all $t>T_0$, and therefore $X^{(1)}(t)>X^{(0)}(t)$ for all $t>T_0$. 
	
	Next, assume that for some $n$ we have $X^{(n)} >  X^{(n-1)}$ and 
$Y^{(n)} < Y^{(n-1)}< 0$. 
	Then, by using first (\iteraX), next (\iteraF) and
(\iteraXnull), then the induction hypotheses $X^{(n)} >  X^{(n-1)}$
and $Y^{(n)} < Y^{(n-1)}$, noting the negative sign in front of the 
integral, then once again the induction 
hypothesis $X^{(n)} > X^{(n-1)}$  but now together with 
$Y^{(n-1)} < 0$ and the negative sign in front of the integral, 
we find for all $t>T_0$ that
$$
\eqalignno{
\hskip-0.5truecm 
	X^{(n+1)}(t) -  X^{(n)}(t) 
& = 
	F_\eps\big(X^{(n)},Y^{(n)}\big)(t) 
	- F_\eps\big(X^{(n-1)},Y^{(n-1)}\big) (t)
&\cr
& = 	
	(1 - \eps)\big(X^{(n)} - X^{(n-1)}\big)(t)
&\cr
& \qquad - \eps  \int_{T_0}^t\int_{T_0}^s 
 \left(Y^{(n)}(s^\pr)e^{2X^{(n)}(s^\pr)} - 
	Y^{(n-1)}(s^\pr)e^{2X^{(n-1)}(s^\pr)}	\right) \dss\ds
&\cr
& \geq
	- \eps  \int_{T_0}^t\int_{T_0}^s 
 \left(Y^{(n)}(s^\pr)e^{2X^{(n)}(s^\pr)} - 
	Y^{(n-1)}(s^\pr)e^{2X^{(n-1)}(s^\pr)}\right) \dss\ds
&\cr
& \geq
	- \eps  \int_{T_0}^t\int_{T_0}^s 
   Y^{(n-1)}(s^\pr) \left( e^{2X^{(n)}(s^\pr)} - e^{2X^{(n-1)}(s^\pr)}\right) 
	\dss\ds
&\cr
& \geq
	0.
&\eqlbl\iteraXmono}
$$
	Hence it follows that $n\mapsto X^{(n)}(t)$ 
is increasing, pointwise for each $t>T_0$.
	Similarly, by using first (\iteraY) and next (\iteraG) and
(\iteraYnull), then the induction hypothesis $Y^{(n)} < Y^{(n-1)}$,
then the induction hypothesis $X^{(n)} >  X^{(n-1)}$ together with 
the negative sign in front of the integral, we find for all $t>T_0$ that
$$
\eqalignno{
\hskip-1truecm 
	Y^{(n+1)}(t) -  Y^{(n)}(t) 
& = 
	G_\eps\big(X^{(n)},Y^{(n)}\big)(t) 
	- G_\eps\big(X^{(n-1)},Y^{(n-1)}\big) (t)
&\cr
& = 	
	(1 - \eps)\big(Y^{(n)} - Y^{(n-1)}\big)(t)
&\cr
& \qquad - \eps  \int_{T_0}^t\int_{T_0}^s {1\over 2}
 \left(e^{2X^{(n)}(s^\pr)} - e^{2X^{(n-1)}(s^\pr)}	\right) \dss\ds
&\cr
& \leq
	- \eps  \int_{T_0}^t\int_{T_0}^s {1\over 2}
	 \left(e^{2X^{(n)}(s^\pr)} - e^{2X^{(n-1)}(s^\pr)}\right) \dss\ds
&\cr
& \leq
	0,
&\eqlbl\iteraYmono}
$$
and it follows that $n\mapsto Y^{(n)}(t)$ is decreasing for each $t>T_0$.\qed
\smallskip

\noindent
{\bf Proposition 3.11:} {\it 
	The joint iteration sequences (\iteraX),
(\iteraY) with initial data (\iteraXnull), (\iteraYnull) 
converge in the limit $n\to\infty$ to asymptotically linear solutions 
of (\NewtonXInew), (\NewtonETAnew) that satisfy (\ACthree) and (\ACfour).}

\smallskip
{\it Proof of Proposition 3.11}:  
	The initial data  $X^{(0)}(t)$ and $Y^{(0)}(t)$ are linear
functions of $t$, with $t>T_0$. 
	We now show  first that a linear upper bound on $X^{(n)}(t)$ 
together with a linear lower bound 
on $Y^{(n)}(t)$ implies corresponding linear bounds on 
$X^{(n+1)}(t)$ and $Y^{(n+1)}(t)$. 
	We then show that these 
bounds converge with $n\to\infty$ to uniform linear bounds for all
$X^{(n)}$ and $Y^{(n)}$. 
	These uniform linear bounds together with the 
monotonicity of the coupled iteration sequences (\iteraX),
(\iteraY) stated in Lemma 3.10 imply that the sequences (\iteraX),
(\iteraY) converge. 
	By inspection of (\iteraX), (\iteraY) 
we see at once that the limit functions are solutions of (\NewtonXInew), 
(\NewtonETAnew) for $t\geq T_0$, with initial data satisfying
the stipulated bounds. 
	Therefore the conclusion holds in particular when the 
initial data are obtained from $\xi(t), \eta(t)$ as $t\to T_0^-$, 
and then the solutions $X(t), Y(t)$ for $t>T_0$ 
coincide with the motion on that trajectory for all $t$. 
	Moreover, the convexity of the trajectories for $t$ large enough
Lemma 3.3, now immediately implies that the trajectories 
are asymptotically straight, with the 
motion on them asymptotically linear, 
satisfying (\ACthree) and (\ACfour), as claimed. 

	It thus remains to prove the uniform linear bounds on 
$X^{(n)}$ and $Y^{(n)}$.
	We begin with the observation that, if for some $n$ the
iterates $X^{(n)}$ and $Y^{(n)}$ satisfy the linear bounds
$$
\eqalignno{
	   X^{(n)}(t)   & < \mu_n \times (t-T_0) + X(T_0),
&\eqlbl\iteraSTARTesta\cr
	0> Y^{(n)}(t)   &>  \nu_n \times (t-T_0)\, ,	
&\eqlbl\iteraSTARTestb}
$$
with some positive constants $\mu_n$ and $\nu_n$, then the
iterates $X^{(n+1)}$ and $Y^{(n+1)}$ satisfy the linear bounds
$$
\eqalignno{
	   X^{(n+1)}(t)   & < \mu_{n+1} \times (t-T_0) + X(T_0),
&\eqlbl\iteraNEXTesta\cr
	0> Y^{(n+1)}(t)   &>  \nu_{n+1} \times (t-T_0)\, ,	
&\eqlbl\iteraNEXTestb}
$$
with 
$$
\eqalignno{
\mu_{n+1}& = \mu_n + \eps\left(\dot{X}(T_0)
		- \delta {\nu_n\over\mu_n^2} -\mu_n\right)
,&\eqlbl\MUnplusone\cr
\nu_{n+1} & = \nu_n + \eps \left(\dot{Y}(T_0) 
		+ \delta {1\over\mu_n} -\nu_n\right).
&\eqlbl\NUnplusone}
$$
	Indeed, by the positivity of $\exp$ and by (\iteraSTARTesta), 
we have
$$
{1\over 2} \int_{T_0}^t\int_{T_0}^s e^{2X^{(n)}(s^\pr)}\dss\ds\, 
< {1\over 2} \int_{T_0}^t\int_{T_0}^\infty e^{2X^{(n)}(s^\pr)}\dss\ds 
< - \delta {1\over\mu_n} (t-T_0),
\eqno\eqlbl\interaINTesta
$$
while by the negativity of $Y^{(n)}$ together with the positivity of
$\exp$, and then by (\iteraSTARTestb), we have
$$
 \int_{T_0}^t\int_{T_0}^s Y^{(n)}(s^\pr)e^{2X^{(n)}(s^\pr)}\dss\ds\, 
>
 \int_{T_0}^t\int_{T_0}^\infty Y^{(n)}(s^\pr)e^{2X^{(n)}(s^\pr)}\dss\ds\, 
> \delta {\nu_n\over\mu_n^2} (t-T_0)\, ,
\eqno\eqlbl\interaINTestb
$$
where
$$
4 \delta =  \exp\big(2X(T_0)\big)\, .
\eqno\eqlbl\deltaDEF
$$
	With these estimates the joint iteration maps (\iteraX),
(\iteraY), with $F_\eps$ and $G_\eps$ given by (\iteraF) and (\iteraG), 
now give (\iteraNEXTesta) and (\iteraNEXTestb)  with
(\MUnplusone) and (\NUnplusone) whenever
(\iteraSTARTesta) and (\iteraSTARTestb) hold.

	Hence, to obtain a linear upper bound on $X(t)$ and
a linear lower bound on $Y(t)$, we need to study the coupled 
recurrence relations (\MUnplusone), (\NUnplusone), starting
with initial data 
$$
\eqalignno{
\mu_{0} & =  \dot{X}(T_0) < 0\, ,
&\eqlbl\initDATAalpha\cr
\nu_{0} & = \dot{Y}(T_0) < 0 \, ,
&\eqlbl\initDATAbeta}
$$
satisfying 
$$
\mu_0^2 +\nu_0^2 =1\, .
\eqno\eqlbl\initDATAcircle
$$
	The last constraint follows from (\COM) and (\COMvalue). 
	The recurrence relations are valid from $n=0$ 
on upward as long as $Y^{(n)}<0$. 
	We need to show that for some legitimate $\mu_0$ and $\nu_0$ 
the recurrence relations converge to limits $\mu_\infty$ and $\nu_\infty$ 
in the desired region of the $\mu,\nu$ plane.

	By inspection we recognize equations (\MUnplusone), (\NUnplusone) 
as the forward Euler approximation to a gradient flow with time step 
$\eps$, defined as follows. 
	We conveniently introduce a new, fictitious ``time'' 
variable $\tau \in \RR^+$ and a $\tau$-dependent point 
$\big(\mu(\tau),\nu(\tau)\big) \in \RR^2$, and we 
let ${\rm Grad}$ denote gradient with respect to $(\mu,\nu)$. 
	We also define the potential 
$$
W(\mu,\nu)= {1\over 2} \Big((\mu - \mu_0)^2+(\nu - \nu_0)^2\Big) 
		- \delta {\nu\over \mu}.
\eqno\eqlbl\ficpot
$$ 
	Then the gradient flow in question is given by
$$
\eqalignno{
	{\dd\over \dtau} (\mu,\nu)(\tau) 
& =  
	- {\rm Grad}\, W\bigl((\mu,\nu)(\tau)\bigr),
&\eqlbl\gradflow\cr
	(\mu,\nu)(0) 
& =  
	(\mu_0,\nu_0),
&\eqlbl\gradflowINITdata}
$$
with initial data $(\mu_0,\nu_0)$ in the set
$$
\SS^1_{-,-} = \SS^1 \cap \RR^2_{-,-}\, ,
\eqno\eqlbl\initZdata
$$
where
$$
\RR^2_{-,-} = \bigl\{ (\mu,\nu)\in\RR^2|\, \mu < 0,\, \nu <0\bigr \}.
\eqno\eqlbl\SWquadrant
$$
	If the gradient flow converges to a stable fixed point, starting at 
the initial datum (\gradflowINITdata), then by choosing $\eps$ small enough 
the iteration (\MUnplusone), (\NUnplusone), starting at
(\initDATAalpha), (\initDATAbeta) will likewise converge to the same 
stable fixed point of (\gradflow). 
	If that fixed point is in $\RR^2_{-,-}$
and the flow from $(\mu_0,\nu_0)$ does not leave $\RR^2_{-,-}$, then the
proposition is proved. 
	It therefore suffices to inspect the gradient flow (\gradflow)
for stable fixed points in  $\RR^2_{-,-}$.
 
	Stable fixed points of the gradient flow (\gradflow) are
critical points of $W$ which locally minimize $W$.
	Clearly, the harmonic oscillator part 
$((\mu -\mu_0)^2+(\nu-\nu_0)^2)/2$ has a unique minimum at
$(\mu_0,\nu_0)$, and an elementary perturbation argument shows
that for each $(\mu_0,\nu_0)$ in the admitted set of initial data
there exists a $\delta_0(\mu_0,\nu_0)>0$ such that, if
$\delta < \delta_0$, then $W(\mu,\nu)$ still has a unique minimum 
at  $(\mu_M,\nu_M)(\delta)$ in the south-western quadrant of 
$\mu,\nu$ space, with $\mu_M > \mu_0$ and $\nu_M < \nu_0$.
	Moreover, the map $\mu_0\mapsto \delta_0$ is strictly
monotonic decreasing. 
	On the other hand, the exponential 
map $X(T_0)\mapsto\delta$ given in (\deltaDEF)
tells us that $\delta \to 0$ rapidly when $X(T_0)\to -\infty$. 
	Also, $\mu_0 \to -1$  as $X(T_0)\to -\infty$. 

	Because of (\xiNULLbound), for $\eta_{\rm in}$ large enough
we have  $X(T_0) \ll -1$, so that we have  $\delta\ll 1$ exponentially
small, given in (\deltaDEF).
	Moreover, we have $(\mu_0,\nu_0)\in \SS^1$ with two negative 
components that satisfy the asymptotic bounds (\xiDOTbound) 
and (\etaDOTbound), so that $(\mu_0,\nu_0)$ is exponentially close 
to the point $(-1,0)$. 
	Therefore, for large negative $X(T_0)$, we surely have $\delta
<\delta_0$. 
	It follows that $W(\mu,\nu)$ then has a unique minimum in the 
south-western quadrant, very close to $\mu_0,\nu_0$ itself. 
	Moreover, along the line $\nu =\nu_0$ the $\nu$ component of the 
gradient flow is given by $\delta/\mu <0$, for $\mu<0$.
	Therefore, the gradient flow (\gradflow) with inital datum 
(\gradflowINITdata) satisfying (\initZdata) remains in  $\RR^2_{-,-}$ and 
converges to $(\mu_M,\nu_M)$.
	The existence proof is complete.\qed
\smallskip

	We have thus shown that for sufficiently large $\eta_{\rm in}>0$
there exists a solution with the correct scattering asymptotics 
(\ACone), (\ACtwo), (\ACthree), (\ACfour). 
	We next reprove our uniqueness statement of Corollary 3.7 by
a different argument that will recur in the sequel.
\smallskip

\noindent
{\bf Theorem 3.12:} {\it The solutions $(\xi(t),\ \eta(t))$ to 
(\NewtonXInew), (\NewtonETAnew) with asymptotic data $\xi_{\rm in}$,
$\eta_{\rm in}$ in (\ACone), (\ACtwo) are unique.} 
\smallskip

{\it Proof of Theorem 3.12}:  
	Let $(\xi_1(t),\ \eta_1(t))$ and $(\xi_2(t),\ \eta_2(t))$   
be two pairs of functions that solve  (\NewtonXInew), (\NewtonETAnew) 
with identical  data (\ACone), (\ACtwo).
	We now define $w_\xi(t) =\xi_1(t)-\xi_2(t)$ and
$w_\eta(t) =\eta_1(t)-\eta_2(t)$ and 
set ${\bf u} = (w_\xi,\dot{w}_\xi,w_\eta,\dot{w}_\eta)^T$.
	Note that 
$$
\lim_{t\to -\infty} {\bf u}(t)	= 0.
\eqno\eqlbl\ASYminusinf
$$
	Next, since $w_\xi$ and $w_\eta$ satisfy the differential equations
$$
\eqalignno{
&
	{\dd^2 w_\xi\over \dt^2} 
= 
	- \eta_1 e^{2\xi_1} + \eta_2 e^{2\xi_2}\, ,
&\eqlbl\wXIeq\cr
&
	{\dd^2 w_\eta\over \dt^2}  
= 
	-{1\over 2}  e^{2\xi_1} +{1\over 2}  e^{2\xi_2}\, ,
&\eqlbl\wETAeq\cr}
$$
by the mean-value theorem there exists a  
$\phi(t)\in
\Big(\min\big(\xi_1(t),\xi_2(t)\big),\max\big(\xi_1(t),\xi_2(t)\big)\Big)$
such that we can rewrite the ODE's for $w_\xi$ and $w_\eta$ as 
$$
\eqalignno{
&
	{\dd^2 w_\xi\over \dt^2} 
= 
	- w_\eta e^{2\xi_1} -  2 w_\xi\eta_2 e^{2\phi}\, ,
&\eqlbl\wXIeqLIN\cr
&
	{\dd^2 w_\eta\over \dt^2}  
= 
	-  w_\xi e^{2\phi}\, .
&\eqlbl\wETAeqLIN\cr}
$$
	We remark that (\wXIeqLIN) and (\wETAeqLIN) are linear 
equations for $w_\xi$ and $w_\eta$.
	We now rewrite (\wXIeqLIN) and (\wETAeqLIN) into the first
order system $\dot{\bf u} = {\bf A}{\bf u}$, where 
$$
{\bf A} = \pmatrix{0		    &1&0	 &0\cr
		  -2\eta_2 e^{2\phi}&0&-e^{2\xi_1}&0\cr
		   0		    &0&0	 &1\cr
		   - e^{2\phi}      &0&0	 &0}
$$
is the coefficient matrix.
	Notice that $\det {\bf A} = -\exp(2\phi+2\xi) < 0$; whence
${\bf  A}$ is invertible.
	More specifically, the characteristic polynomial of ${\bf A}$
is readily found to be
$$
P(\lambda)=\lambda^4 + 2{\eta_2}e^{2\phi}\lambda^2 - e^{2\phi+2\xi_1}.
$$
	Solving for the roots of $\lambda^2$ we find two real values
$$
\lambda^2 =  
\Big(-{\eta_2} \pm \sqrt{{\eta_2}^2 + e^{2\xi_1}}\Big)e^{\phi} ,
$$
one positive, the other negative.
	Hence, there are 2 real and 2 purely imaginary eigenvalues $\lambda$
of ${\bf A}$. 
	Now, in view of (\ASYminusinf)
the purely imaginary roots do not contribute to the solutions with
our scattering data.
	Next, the real roots are
$$
\lambda^R_\pm = 
\pm\Big(-{\eta_2} + \sqrt{{\eta_2}^2 + e^{2\xi_1}}\Big)^{1\over 2}e^{\phi/2},
$$
one negative, the other positive for all $t\in \RR$.
	Thus, $\phi(t)\sim -|t|$ for  $t\to - \infty$, 
by letting $t\to - \infty$ we see that the real roots converge to 0
exponentially fast.
	Hence the nontrivial orbits of $\dot{\bf u}={\bf A}{\bf u}$ coming 
from the real roots converge to some ${\bf u}^\sharp\neq {\bf 0}$ ouside 
some ball in $\RR^4$, centered at the origin. 
	Therefore, the only vector solution compatible with the asymptotic 
conditions (\ASYminusinf) is ${\bf u} \equiv {\bf 0}$, viz.
$w_\xi(t) \equiv 0 \equiv w_\eta(t)$.
	Uniqueness is proved.\qed
\smallskip
\vfill\eject

	We remark that Theorem 3.12, like Corollary 3.7, claims
uniqueness not only for the scattering solutions for which there 
exists a $\Theta \in (-\pi,-\pi/2)$. 
	We now return to those scattering solutions and 
show that there exist scattering solutions for the whole
range of deflection angles $\Theta\in (-\pi,-\pi/2)$. 
\smallskip

\noindent
{\bf Theorem 3.13:} {\it For every $\Theta\in (-\pi,-\pi/2)$
there is a choice of parameters $\eta_{\rm in}>0$ and $\xi_{\rm in}$
such that there exists a solution $(\xi(t),\ \eta(t))$ to 
 (\NewtonXInew), (\NewtonETAnew) with scattering data 
(\ACone), (\ACtwo), (\ACthree), (\ACfour).}
\smallskip

{\it Proof of Theorem 3.13}: 
	We argue via continuity. 
\smallskip

\noindent
{\bf Definition 3.14:} {\it We define $S$ to be the set 
$(\xi_{\rm in},\eta_{\rm in}, \Theta)\in\RR^3$ for which there 
exists a joint solution
$\xi=f_\xi(t),\ \eta = f_\eta(t)$ of (\NewtonXInew), (\NewtonETAnew) 
satisfying the asymptotic conditions 
(\ACone), (\ACtwo), (\ACthree), (\ACfour).
}
\smallskip

	Let $\RR^+  = (0,\infty)$ and set
$
W = \RR\times\RR^+\times(-\pi,-\pi/2).
$ 
	We will show that $S$ is relatively open and closed in $W$.
	Clearly, by our existence proof, $S$ is non-empty; thus, $S$
is a connected non-empty set and it will follow   that the projection
of $S$ onto the third component is $(-\pi,-\pi/2)$.
	To show that $S$ is open we will apply the inplicit function 
theorem to our ODEs  (\NewtonXInew), (\NewtonETAnew),
fix $s_0\in S$, and we have a solution $\xi_0(t),\eta_0(t)$
with  scattering data $s_0$.
\smallskip

	To show that $S$ is open, we consider the
linearized part of  $\xi = \xi_0 + \xi_1 + ...$ and
$\eta = \eta_0 +\eta_1 + ...$, with $\xi_1$ and $\eta_1$ small, satisfying
$$
\lim_{|t|\to\infty} \xi_1(t) = 0=  \lim_{|t|\to\infty} \eta_1(t)
\eqno\eqlbl\phipsiASYMPcond
$$
and satisfying the linearized equations of motion
$$
\eqalignno{
&
	{\dd^2\xi_1\over \dt^2} 
= 
	- \eta_1 e^{2\xi_0} -2 \xi_1\eta e^{2\xi_0}\, ,
&\eqlbl\NewtonXInewLIN\cr
&
	{\dd^2\eta_1\over \dt^2}  
= 
	- \xi_1 e^{2\xi_0}\, .
&\eqlbl\NewtonETAnewLIN\cr}
$$
	Rewriting these second order equations as
a first order system for ${\bf v}^T = (\xi_1,\eta_1)$, 
we are led to $\dot{\bf v}={\bf M}{\bf v}$, with coefficient matrix
$$
{\bf M} = \pmatrix{0		        &1&0	      &0\cr
		   -2 \eta_0 e^{2\xi_0} &0&-e^{2\xi_0}&0\cr
		   0		        &0&0	      &1\cr
		   - e^{2\xi_0}	        &0&0	      &0}
$$
and with ${\bf v}(t)\to {\bf 0}$ as $|t|\to\infty$.
	Clearly,  similar to the proof of Theorem 3.12, we have 
$\det{\bf M} = -\exp(4\xi_0) < 0$, and the characteristic polynomial is
$$
P(\lambda) = \lambda^4 + 2{\eta_0}e^{2\xi_0}\lambda^2 - e^{4\xi_0},
$$
with 2 real and 2 purely imaginary eigenvalues $\lambda$
of ${\bf M}$, for all $t\in \RR$. 
	Thus, by the condition that ${\bf v}(t)\to {\bf 0}$
for $|t|\to \infty$, we conclude that 
${\bf u}(t) ={\bf 0}$ identically.
	Therefore, the implicit function theorem applies and we may
conclude that there is a neighborhood about $s_0$ in $W$ for which one
finds solutions to   (\NewtonXInew), (\NewtonETAnew), 
satisfying the asymptotic conditions 
(\ACone), (\ACtwo), (\ACthree), (\ACfour).
	Hence, $S$ is an open set.

	To show that $S$ is relatively closed, consider a sequence
$s_n\in S$ such that $s_n\to s_*\in W$.
	We have 
$s_n =(\xi_{{\rm in},n},\eta_{{\rm in},n},\Theta_n)$
and 
$s^* =(\xi_{{\rm in}}^*,\eta_{{\rm in}}^*,\Theta^*)$.
	Note that we have solutions of
$$
\eqalignno{
&
	{\dd^2\xi_n\over \dt^2} 
= 
	- \eta_n e^{2\xi_n}\, ,
&\eqlbl\NewtonXIn\cr
&
	{\dd^2\eta_n\over \dt^2}  
= 
	- {1\over 2} e^{2\xi_n}\, ,
&\eqlbl\NewtonETAn\cr}
$$
satisfying the scattering data for $s_n$, by the very Definition 3.14
of $S$.

	Because $s_n$ belongs to a bounded set with
compact closure in $W$, by (\xiINT) and (\etaINT) the asymptotic behavior 
of $(\xi_n(t),\eta_n(t))$ in (\ACone), (\ACtwo) 
is uniform, and independent of the solution $(\xi_n,\eta_n)$. 
Similarly we have uniformity in (\ACthree) and (\ACfour)
	That means, the error term is uniform in $n$ if $s_n$ 
remains in a set with compact closure in $W$. Similarly,
by differentiating (3.22) and (3.23) once and
using (\COM) and the uniformity in $(\eta_n(t), \xi_n(t))$ we may 
conclude the same uniformity for the derivatives.
	This allows us to conclude compactness at ``infinity.''  

	First, we conclude that
$$
\sup_{t,n} (\dot\xi_n^2 +\dot\eta_n^2)^{1/2} \leq c.
\eqno\eqlbl\speedlimit
$$
	To see that (\speedlimit) holds, indeed, recall that
$\dot\eta_n$ is strictly monotonic decreasing, by (\NewtonETAnew).
	Since by hypothesis, 
$\lim_{t\to\infty}\dot\eta_n(t) = \sin\Theta_n$ and 
$\lim_{t\to -\infty}\dot\eta_n(t) = 0$, we have that 
$|\dot\eta_n|\leq |\sin\Theta_n|$, but also $\dot\eta_n <0$  and
therefore $\eta_n$ strictly monotonic decreasing. 
	Furthermore, as long as $\eta \geq 0$, we have that 
$\dot\xi_n$ is strictly monotonic decreasing,  by (\NewtonXInew),
and when $\eta_n =0$ at $t=T_0$, we have $\dot\xi_n^2 + \dot\eta_n^2 =1$,
by (\COM) and (\COMvalue). 
	Thus, since also 
$\lim_{t\to -\infty}\dot\xi_n(t) = 1$, we conclude that 
$|\dot\xi_n|\leq 1$  for $t\in (-\infty,T_0]$.
	On the other hand, for $t >T_0$ we have $\eta_n<0$ by the
strict monotonic decrease of $\eta_n$, and thus by (\NewtonXInew) 
we now have that $\dot\xi_n$ is strictly monotonic increasing for
$t> T_0$.
	But then, since 
$\lim_{t\to\infty}\dot\xi_n(t) = \cos\Theta_n$, we conclude that
$|\dot\xi_n|\leq 1$  for $t\in (T_0,\infty)$ as well. 
	Thus, (\speedlimit) is established.

	Next we show that  there is a point $t_n$, with $|t_n|\leq c^\prime$
independent of of $n$, and some $C$ independent of $n$, 
such that 
$$
	|\xi_n(t_n)| +|\eta_n(t_n)| \leq C.
\eqno\eqlbl\pointNORMbound
$$
	Thus, pick 
$t_n = \ln( 2\sqrt{\eta_{{\rm in},n}}) - \xi_{{\rm in},n}$.
	Then, by (\XIupperBOUND), we have
$\xi_n(t_n) \leq \xi_{{\rm in},n} + t_n$.
	We proved in Lemma 3.5 (see also the proof of Lemma 3.9) 
that for $t< T_{1/2}$ we have $\xi(t) > \xi^{(0)}(t)$, with $\xi^{(0)}$ 
given in  (\XInull), 
and this thus holds for any $\xi_n$ with a corresponding $T_{1/2,n}$.
	Thus, since $t_n < T_{1/2,n}$, by (\THALFlowerBOUND), we have
$$
	\xi_n(t_n) 
>  
	-\ln\cosh \big(t_n + \xi_{{\rm in},n} 
	- \ln(2/\sqrt{\eta_{{\rm in},n}})\big)
	- \ln \sqrt{\eta_{{\rm in},n}}\, ,
\eqno\eqlbl\XInBOUND
$$
and the bounds for $\xi_n$ are established. 
	Since $s_n$ belongs to a set with compact closure, it follows that
there exists a $c^\prime$ independent of $n$ such that $|t_n|< c^\prime$. 
	
	Next, we know that $\eta_n$ is a decreasing function, bounded above
by $\eta_n < \eta_{\rm in}$.
	By Lemma 3.8, since $T_{1/2,n}> t_n$, we see that 
$\eta_n(t_n) \geq 	\eta_n(T_{1/2,n})$.
	Thus, $|\eta_n(t_n) |$ is bounded above independent of $n$, too,
and this finishes the proof of (\pointNORMbound).

	Next, using (\speedlimit) and  (\pointNORMbound), we conclude that
$\|(\xi_n,\eta_n) \|_{_{L^\infty(I)}}\leq C(I)$, where $I$ is any
bounded sub-interval of $\RR$.
	Thus, by using (\speedlimit) and the Ascoli theorem we
conclude that $(\xi_n,\eta_n)$ converges uniformly on bounded
sub-intervals of $\RR$ to continuous functions $(\xi^*,\eta^*)$. 
	Using now (\NewtonXIn), (\NewtonETAn), this uniform
convergence now  implies that the second derivatives 
$(\ddot\xi_n,\ddot\eta_n)$ are uniformly bounded on compact 
sub-intervals of $\RR$.
	Since we also have (\speedlimit), by Ascoli's theorem
again, the first derivatives $(\dot\xi_n,\dot\eta_n)$ 
converge uniformly to $(\dot\xi^*,\dot\eta^*)$ 
bounded on compact sub-intervals of $\RR$.
	Therefore, in the sense of distributions, 
$$
\eqalignno{
&
	{\dd^2\xi_n\over \dt^2} 
=
	- \eta^* e^{2\xi^*}\, ,
&\eqlbl\NewtonXInn\cr
&
	{\dd^2\eta_n\over \dt^2}  
=
	- {1\over 2} e^{2\xi^*}\, .
&\eqlbl\NewtonETAnn\cr}
$$

	Next, we readily establish that 
$\lim_{t\to\infty} t^{-1}\xi^*(t) = \cos\Theta^*$, 
that $\lim_{t\to\infty} t^{-1}\eta^*(t) =\sin\Theta^*$,
and also that
$\lim_{t\to\infty} \xi^*(t)-t = \xi_{{\rm in}}^*$,
and $\lim_{t\to\infty} \eta^*(t) = \eta_{{\rm in}}^*$.
	Thus $(\xi^*(t),\eta^*(t))$ satisfies the asymptotic 
conditions (\ACone), (\ACtwo), (\ACthree), (\ACfour); 
hence, $(\xi^*,\eta^*)$ is a solution, and therefore
$S$ is open and relatively closed in $W$.

	Since $W$ is connected and $S\neq\emptyset$, we conclude that
$S$ is a connected set in $W$.
	To finish the proof, we need to show that the projection of 
$S$ onto the third component of $W$ is indeed the full interval
$(-\pi, -\pi/2)$. 
	Since $S$ is connected and open, and since the projection map
is continuous and open, the projection of $S$ into $(-\pi, -\pi/2)$
is an interval, say $(\vartheta_1,\vartheta_2)$, with
$-\pi < \vartheta_1$ and  $\vartheta_2 < - \pi/2$.
	Thus, for instance, as $\Theta_j\to\vartheta_1$,
either $\eta_{{\rm in},j}\to 0$ or $\eta_{{\rm in},j}\to \infty$.
	Let $\eta_j(t_j) \to 0$.
	Assuming that  $\eta_{{\rm in},j}\to \infty$ as 
$\Theta_j\to\vartheta_1$, from (\xiINT)  we conclude that
$\xi_j(t_j)\to -\infty$, which now contradicts the condition
that $\Theta_j\to \vartheta_1\in (-\pi, -\pi/2)$.
	Assuming that  $\eta_{{\rm in},j}\to 0$ as 
$\Theta_j\to\vartheta_1$, we again arrive at the 
contradiction by Lemma 3.4. 
	The other cases are $\xi_{{\rm in},j} \to\pm \infty$
for fixed $\eta_{\rm in}$.
	Assume first that $\xi_{{\rm in},j} \to - \infty$.
	Then by (\THALFlowerBOUND) we see that $T_{1/2}\to +\infty$
for fixed $\eta_{\rm in}$, which means that $\eta(t)> \eta_{\rm in}/2$
for all $t\in \RR$, which is impossible.
	Finally, assume that $\xi_{{\rm in},j} \to + \infty$,
for fixed $\eta_{\rm in}$. 
	Then, since $\xi^{(0)}$ is a subsolution for $\xi$, 
we have that 
$$
	\eta(t) 
< 
	\eta_{\rm in} -{1\over 2} \int_{-\infty}^t\int_{-\infty}^s
	e^{2\xi^{(0)}(\ts)}\dd \ts\dd s,
\eqno\eqlbl\etaUPbound
$$
for all $t$. 
	Using (\XInull), we obtain
$$
	\eta(t) 
< 
	\eta_{\rm in} -e^{2\xi_{\rm in}} F(t),
\eqno\eqlbl\etaUPboundb
$$
where $F(t)$ is a monotonically increasing, positive function,
and $F(t)\to 0$ exponentially fast as $t\to - \infty$.
	Next, let $T_{0,j}$ be defined by $\eta(T_{0,j})=0$. 
	Clearly, we now conclude from (\etaUPboundb) and the
properties of $F$ that $T_{0,j}\to -\infty$ as $\xi_{{\rm in},j}\to +\infty$.
	But then, we conclude that $-1\leq \dot\eta(t) < 0.5\sin\vartheta_1$
for all $t> T_{0,j}$, 
with $T_{0,j}\to -\infty$ as $\xi_{{\rm in},j}\to +\infty$, in 
contradiction to (\ACtwo).	
 
	This concludes our proof of Theorem 3.13. \qed 
\smallskip
	The proof of Theorem 3.1 is complete. \qed 

\noindent
{\bf Acknowledgement:} S. Chanillo was supported by the NSF through
Grant DMS-9970359.
\vfill\eject

\biblio

\vfill\eject
\epsfxsize=8.5cm
\centerline{\epsffile{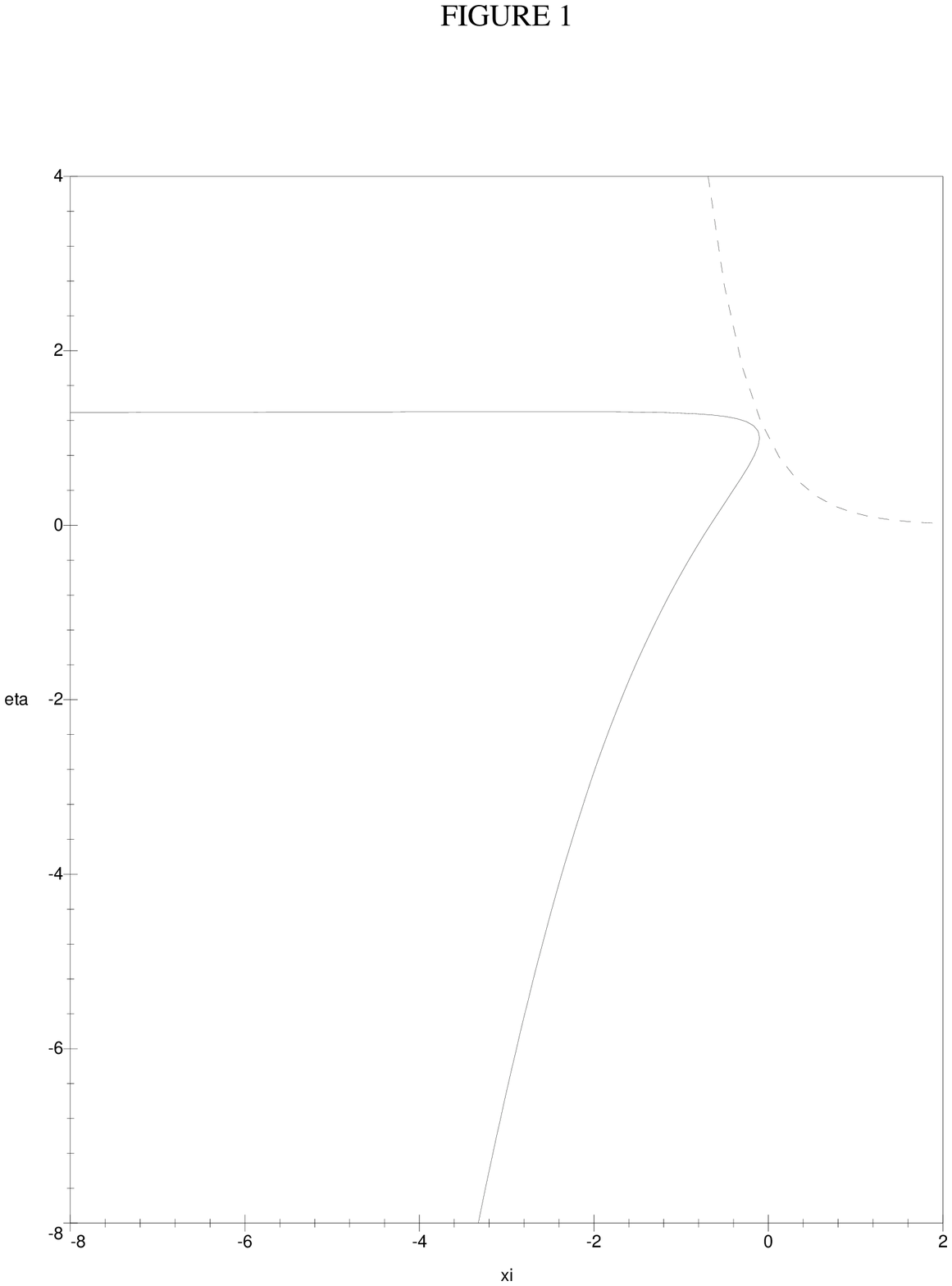}}

\noindent
Fig.1: A regular scattering trajectory (solid curve) with relevant
scattering data. For convenience, the locus of singular points 
(dashed curve) is displayed as well.
\bigskip
\bigskip
\vfill\eject

\medskip
\epsfxsize=8.5cm
\centerline{\epsffile{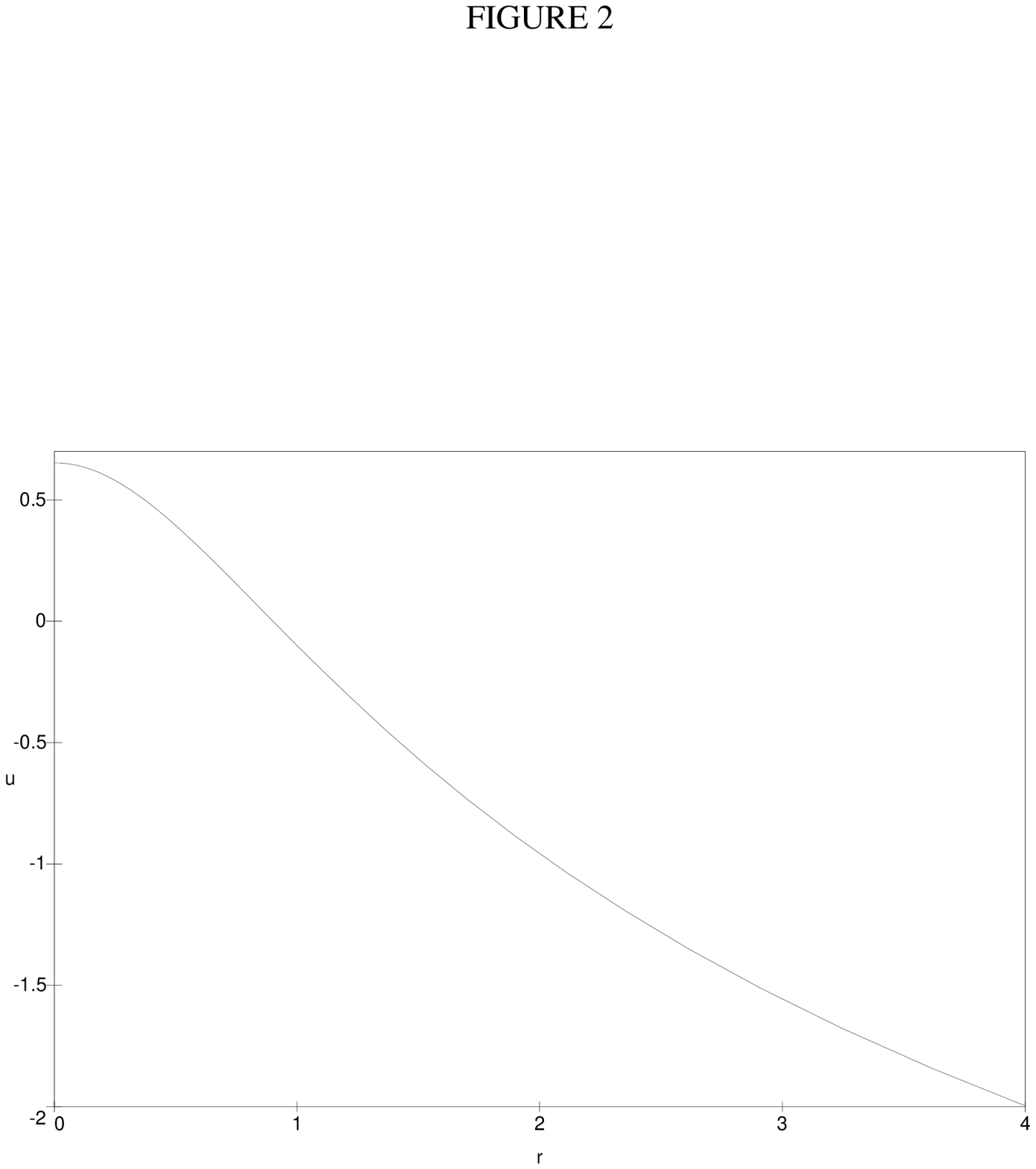}}

\noindent
Fig.2: The solution $u$ as function of $r$ obtained from the motion on
the scattering trajectory of Fig.1.
\bigskip
\bigskip

\vfill\eject

\medskip
\epsfxsize=8.5cm
\centerline{\epsffile{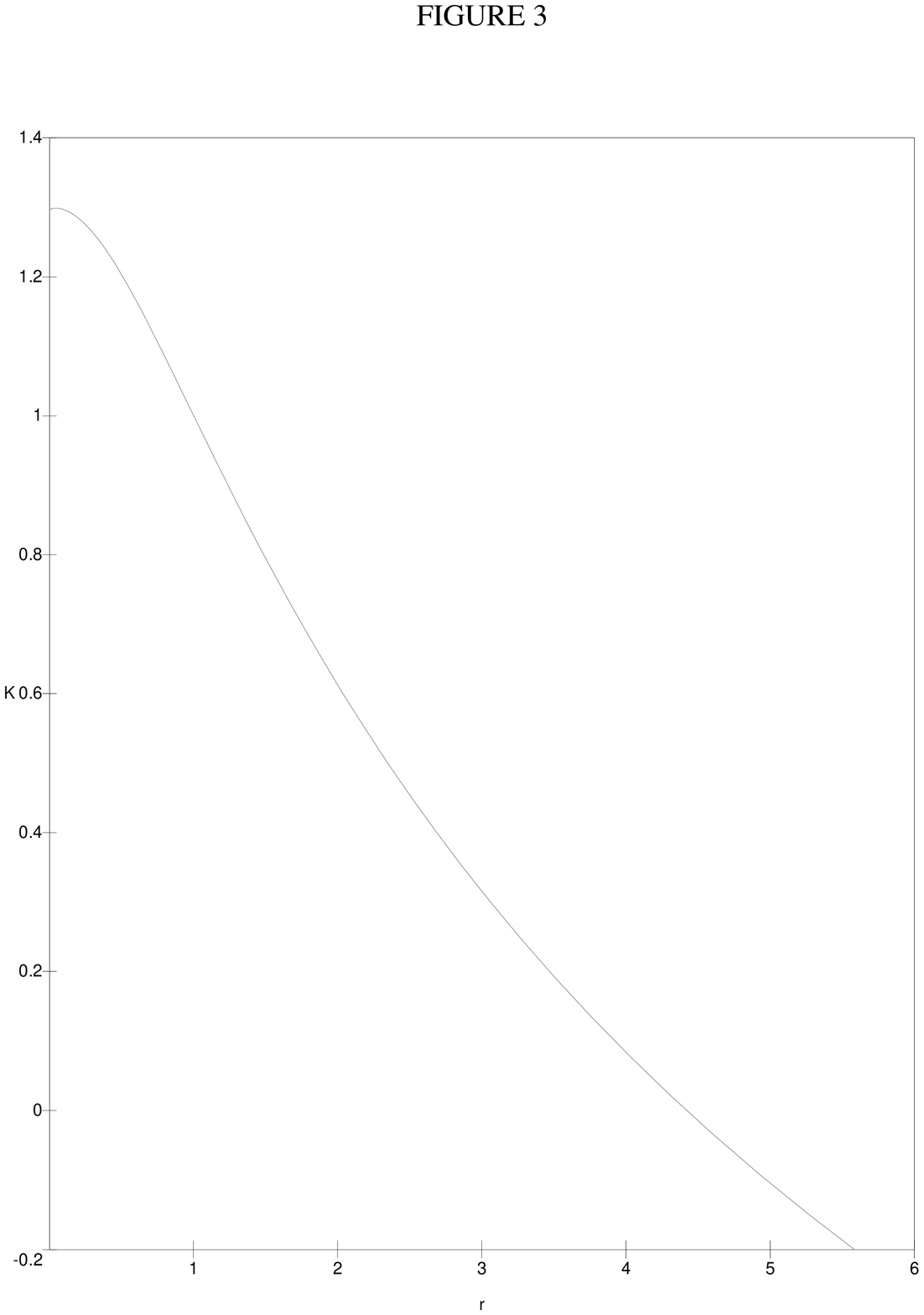}}

\noindent
Fig.3: The Gauss curvature $K$ as function of $r$ obtained from the motion 
on the scattering trajectory of Fig.1.
\bigskip
\bigskip

\bye